\documentclass[11pt]{amsart}
\usepackage{amsmath,amsthm, amscd, amssymb, amsfonts}
\usepackage[all]{xy}

\newcommand{\op}{{\rm op}}

\newcommand\boxe{\begin{tabular}{|p{0,1cm}|}
\hline \\ \hline \end{tabular}}

\newcommand\boxee{\begin{tabular}{|p{0,8cm}|}
\hline \\ \\ \hline \end{tabular}}

\newcommand{\lyzu}{\xi}
\newcommand{\lyzd}{\eta}

\newcommand{\lyzuu}{\operatorname{in}_1}
\newcommand{\lyzdd}{\operatorname{in}_2}

\newcommand\Gpd{\operatorname{Gpd}}
\newcommand\Quiv{\operatorname{Quiv}}
\newcommand\supp{\operatorname{supp}}
\newcommand\ev{\operatorname{ev}}
\newcommand\coev{\operatorname{coev}}

\DeclareMathOperator*{\Tim}{\times}
\newcommand{\ftimes}{\sideset{_\tgt}{_\sou}\Tim}

\newcommand{\prode}{{\boxtimes}}
\newcommand{\gstr}{\Gb_{\Ac}}
\newcommand{\hstr}{\Hb_{\Ac}}
\newcommand{\gstro}{\Gb_{\Ac^{\op}}}
\newcommand{\ggstr}{\G_{\Ac}}
\newcommand{\sou}{\mathfrak s}
\newcommand{\tgt}{\mathfrak e}
\newcommand{\conexo}{\approx}

\newcommand{\R}{{\mathcal R}}

\newcommand{\W}{{\mathcal W}}
\newcommand{\Uc}{{\mathcal U}}

\newcommand{\q}{\mathbf{q}}
\newcommand{\ub}{\mathbf{u}}
\newcommand{\w}{\mathbf{w}}

\newcommand{\ku}{\Bbbk}
\newcommand{\K}{{\mathcal K}}
\newcommand{\Ac}{{\mathcal A}}
\newcommand{\Bc}{{\mathcal B}}
\newcommand{\N}{{\mathcal N}}
\newcommand{\uno}{{\bf 1}}
\newcommand{\path}{\operatorname{Path}}
\newcommand{\bundle}{\operatorname{bundle}}
\newcommand{\Rep}{\operatorname{Rep}}

\newcommand{\f}{\varphi}

\newcommand{\x}{{\overline x}}
\newcommand{\y}{{\overline y}}
\newcommand{\z}{{\overline z}}
\newcommand{\osigma}{{\overline \sigma}}
\newcommand{\psigma}{\sigma^{\rm P}}

\newcommand{\X}{{\mathcal X}}
\newcommand{\G}{{\mathcal G}}

\newcommand{\F}{{\mathcal F}}
\newcommand{\C}{{\mathcal C}}
\newcommand{\D}{{\mathcal D}}
\newcommand{\Ec}{{\mathcal E}}

\newcommand{\m}{\mathcal{M}}
\newcommand{\n}{\mathcal{N}}
\newcommand{\Lc}{\mathcal{L}}

\newcommand{\T}{{\mathcal T}}
\newcommand{\Hc}{{\mathcal H}}
\newcommand{\Vc}{{\mathcal V}}
\newcommand{\Pc}{{\mathcal P}}

\newcommand{\Hb}{{\mathbb H}}

\newcommand{\Gb}{{\mathbb G}}

\newcommand{\Ag}{{\mathfrak A}}
\newcommand{\Bg}{{\mathfrak B}}
\newcommand{\Pg}{{\mathfrak P}}
\newcommand{\Hg}{{\mathfrak H}}
\newcommand{\Vg}{{\mathfrak V}}

\newcommand{\aut}{\mathfrak{aut}\,}

\newcommand{\Lin}{\operatorname{Lin}}

\newcommand{\Res}{\operatorname{Res}}

\newcommand{\Inn}{\operatorname{Inn}}

\newcommand\card{\operatorname{card}}

\newcommand\Hom{\operatorname{Hom}}

\newcommand\Ker{\operatorname{Ker}}
\newcommand\Imm{\operatorname{Im}}

\newcommand{\trid}{\,{\triangleright}\,}
\newcommand{\fde}{\,{\rightharpoonup}\,}
\newcommand{\fiz}{\,{\leftharpoonup}\,}

\newcommand{\fd}{\,{\rightharpoondown}\,}
\newcommand{\fz}{\,{\leftharpoondown}\,}

\newcommand{\fdd}{\,{\hookrightarrow}\,}
\newcommand{\fzz}{\,{\hookleftarrow}\,}

\newcommand{\fddd}{\,{\hookrightarrow}\,}
\newcommand{\fzzz}{\,{\hookleftarrow}\,}

\newcommand\Mpp{{\,}_{\Pc} \hspace{-0.025cm}\mathcal M_{\Pc}}
\newcommand\mpp{{\,}_{\Pc} \hspace{-0.025cm}\mathfrak m_{\Pc}}

\numberwithin{equation}{section}\theoremstyle{plain}

\newtheorem{step}{Step}
\newtheorem{theorem}{Theorem}[section]
\newtheorem{lema}[theorem]{Lemma}
\newtheorem{cor}[theorem]{Corollary}

\newtheorem{proposition}[theorem]{Proposition}

\theoremstyle{definition}
\newtheorem{definition}[theorem]{Definition}
\newtheorem{exa}[theorem]{Example}
\newtheorem{exas}[theorem]{Examples}

\theoremstyle{remark}
\newtheorem{obs}[theorem]{Remark}

\newcommand\id{\operatorname{id}}
\newcommand\idd{\mathbf{id}}

\def\pf{\begin{proof}}
\def\epf{\end{proof}}

\theoremstyle{remark}

\begin{document}

\renewcommand{\baselinestretch}{1.2}
\thispagestyle{empty}
\title[On the quiver-theoretical QYBE]{On the quiver-theoretical  quantum Yang-Baxter equation}
\author[Nicol\'as Andruskiewitsch]{Nicol\'as Andruskiewitsch}
\address{Facultad de Matem\'atica, Astronom\'\i a y F\'\i sica
\newline \indent
Universidad Nacional de C\'ordoba
\newline
\indent CIEM -- CONICET
\newline
\indent (5000) Ciudad Universitaria, C\'ordoba, Argentina}
\email{andrus@mate.uncor.edu, \quad \emph{URL:}\/
http://www.mate.uncor.edu/andrus}
\thanks{Partially supported by CONICET,
Agencia C\'ordoba Ciencia, ANPCyT    and Secyt (UNC)}
\date{\today}
\begin{abstract}  Quivers over a fixed base set
form a monoidal category  with tensor product given by pullback.
The quantum Yang-Baxter equation, or more properly the braid
equation, is investigated in this setting. A solution of the braid
equation in this category  is called a ``solution" for short.
Results of Etingof-Schedler-Soloviev, Lu-Yan-Zhu and Takeuchi on
the set-theoretical  quantum Yang-Baxter equation are generalized
to the context of quivers, with groupoids playing the r\^ole of
groups. The notion of ``braided groupoid" is introduced. Braided
groupoids are solutions and are characterized in terms of
bijective 1-cocycles. The structure groupoid of a non-degenerate
solution is defined; it is shown that it is braided groupoid. The
reduced structure groupoid of a non-degenerate solution is also
defined. Non-degenerate solutions are classified in terms of
representations of matched pairs of groupoids. By linearization we
construct star-triangular face models  and realize them
as modules  over
quasitriangular quantum groupoids introduced in recent papers by
M. Aguiar, S. Natale and the author. \end{abstract}

\maketitle

\section*{Introduction}

The problem of classification of set-theoretical solutions to the
quantum Yang-Baxter equation (QYBE) was raised by Drinfeld
\cite{Dr}, who also mentioned the example of a subset of a group
stable under conjugation. The question was considered by
Etingof-Schedler-Soloviev \cite{ess, s} and Lu-Yan-Zhu
\cite{lyz1}. In these works, an abstract characterization of
solutions to the QYBE in group-theoretical terms was given.
Indecomposable solutions with an underlying set with a prime number of elements
were classified in \cite{egs}. Later, Takeuchi simplified some
arguments and provided a unified presentation in terms of braided
groups \cite{tak}.

\bigbreak
It is natural to extend Drinfeld's question and ask for solutions
to the QYBE in an arbitrary monoidal category.
In the linear or in the set-theoretical settings, or more
generally in a symmetric monoidal category, the QYBE is equivalent
to the braid equation. Actually, the main goal of the above
mentioned articles is the study of the braid equation-- and
information on the QYBE is obtained as a bonus. Such an
equivalence does not exist in arbitrary, not symmetric, monoidal categories.

\bigbreak In this paper, the question of solutions to the braid
equation is addressed in the special case of the category of
quivers over a fixed set $\Pc$ of vertices; we refer to this as
the quiver-theoretical braid equation. Results of
Etingof-Schedler-Soloviev \cite{ess, s} and Lu-Yan-Zhu \cite{lyz1}
on the set-theoretical  quantum Yang-Baxter equation are
generalized to the context of quivers, with groupoids playing the
r\^ole of groups. Our approach, inspired by the presentation of
Takeuchi \cite{tak}, is different from the original papers
\cite{ess, lyz1, s}; it is systematically based on the notion of
braided groupoid, whereas in those papers much emphasis was put on
the derived structure groupoid.

\bigbreak In Section 1, basic notions on quivers and groupoids are
reminded. The apparently new construction of the free groupoid
generated by a quiver is given. This is needed in Section 2 to
define the structure groupoid of a non-degenerate solution. In
Section 2, braided groupoids are investigated. The main results on
characterization of non-degenerate solutions are stated and proved
in Section 3. Our main result, Theorem \ref{finito},  gives a
classification of non-degenerate solutions via representations of
suitable matched pairs of groupoids. This result generalizes
\cite[Th. 2.7]{s} but the formulation is new even in the
set-theoretical case. The concept of representations of  matched
pairs of groupoids was introduced in \cite{AA}. In this paper,
``positive" universal $R$-matrices for the quantum groupoids
constructed in \cite{AN} are presented. In Section 4, devoted to
linearization, we construct star-triangular face models and
realize some of them over quasitriangular quantum groupoids from
\cite{AN, AA}.

\bigbreak We close this introduction with three remarks.

\bigbreak
The category of quivers
over $\Pc$, with fiber product as tensor product, does not appear
to be symmetric in a reasonable way. However there are some
substitutes of the symmetry and a formulation of the
quiver-theoretical quantum Yang-Baxter equation is still possible
\eqref{qtQYBE}, and there
is still an equivalence between solutions to the braid equation
and solutions to \eqref{qtQYBE}.

\bigbreak
The analogy of some of the results in the set-theoretical and
quiver-theoretical settings suggests that these might be
particular cases of a general description of solutions to the
braid equation in arbitrary monoidal categories with extra
hypothesis (\emph{e. g.} existence of equalizers and direct
products).

\bigbreak After release of the first version of this paper, it was
pointed out to the author that a solution of the set-theoretical
Yang-Baxter equation is essentially the same thing as a birack;
see \cite{ces, fjk, frs, sw, w} and references therein.

\subsection*{Acknowledgements} This work is
a continuation of \cite{AN, AA}. Aprovecho esta ocasi\'on para
agradecer a Sonia y Marcelo por compartir conmigo su entusiasmo e
intuici\'on. I also thank Professor M. Takeuchi for sending me a
copy of \cite{tak}; the strategy in the present paper owes a lot
to it.

\section{Quivers and groupoids}

\bigbreak
\subsection{Quivers}

\

\bigbreak Let  $(\Ac, \Pc, \sou, \tgt)$ be a quiver; thus $\Ac$
and $\Pc$ are sets, with $\Pc$ non-empty, and $\sou, \tgt: \Ac
\rightrightarrows \Pc$ are functions. An element $a$ of $\Ac$ is
an ``arrow" from its source $\sou(a)$ to its end $\tgt(a)$. We
shall fix $\Pc$ and say that ``$\Ac$ is a quiver" or ``$\Ac$ is a
quiver over $\Pc$".
Quivers are also called ``oriented graphs".
We denote by $\Ac(P, Q)$ the set of arrows
from $P$ to $Q$, if $P, Q\in \Pc$; and $\Ac(P) = \Ac(P, P)$.
Morphisms of quivers over $\Pc$ are defined in the usual way; they
should be the identity on $\Pc$.

\bigbreak A quiver $\Bc$ differs from $\Ac$ in the orientation if
it has the same arrows as $\Ac$ but with different ${\mathfrak
s}$, $\tgt$ for some arrows. For instance, the
\emph{opposite quiver} is  $\Ac^{\op} = \Ac \times \{-1\}$; if
$x\in \Ac$ then $x^{-1} := (x,-1)$ has $\mathfrak s(x^{-1}) =
\mathfrak e (x)$, $\mathfrak e(x^{-1}) = \mathfrak s(x)$. Also, by
abuse of notation, we set $(\Ac^{\op})^{\op} = \Ac$, and
$(x^{-1})^{-1} = x$ for $x\in \Ac$.

\bigbreak If $\Ac$ and $\Bc$ are quivers over $\Pc$ then we can
form the disjoint union $\Ac \coprod \Bc$ that is again a quiver over
$\Pc$. The \emph{double} of $\Ac$ is $\D \Ac  := \Ac \coprod
\Ac^{\op}$; it does not depend on the orientation of $\Ac$.
The quiver $\D\Ac$ is occasionally denoted by $\overline{\Ac}$ in the literature.

\bigbreak Let $n$ be a positive integer. A \emph{path of length
$n$} in $\Ac$ is a sequence $w = (x_{1}, \dots, x_{n})$ of elements in
$\Ac$ such that $\sou (x_{i+1}) = \tgt (x_{i})$, $1\le i < n$; we
shall denote it by $w = x_{1} x_{2} \dots x_{n}$. A \emph{path
of length 0} is a symbol $\id P$, $P\in \Pc$. The set of all
paths of length $n$ in $\Ac$ is a quiver $\path_n(\Ac)$
with $\sou (w) = \sou (x_{1})$,
$\tgt (w) = \tgt (x_{n})$ if $w = x_{1} x_{2}
\dots x_{n}$ (if $n>0$), and $\sou (\id P) = \tgt (\id P) =
P$ (if $n=0$). The quiver of all paths in $\Ac$
is $\path(\Ac) = \coprod_{n\ge 0}\path_n(\Ac)$.

\bigbreak The quiver $\Ac$ induces an equivalence relation on
$\Pc$: $P\conexo Q$ if and only if there exists $w\in \path
(\D \Ac )$ with $\sou (w) = P$, $\tgt (w) = Q$,
$P, Q\in \Pc$. Then $\Ac$ is \emph{connected} if $P\conexo Q$ for all
$P ,Q\in \Pc$.

\bigbreak Any map $p: \Lc \to \Pc$ can be considered as a quiver
with $\sou = \tgt = p$; such a quiver shall be called a \emph{loop
bundle}. In this context we shall sometimes use the fiber
notation: $\Lc_P := \Lc(P) = p^{-1}(P)$.

A quiver $\Ac$ gives rise to two loop bundles: $\sou: \Ac \to \Pc$
and $\tgt: \Ac \to \Pc$; we shall denote them by $\Ac^{\sou}$ and
$\Ac^{\tgt}$, respectively. It might be useful to visualize them
as follows:
\begin{align} \label{loop1}
\Ac^{\sou} &= \coprod_{P \in \Pc} \Ac^{\sou}(P), \qquad
\Ac^{\sou}(P) = \{(y, y^{-1}): y\in \Ac(P, Q), \, Q\in \Pc\},
\\ \label{loop2} \Ac^{\tgt} &= \coprod_{Q \in \Pc} \Ac^{\tgt}(Q),
\qquad \Ac^{\tgt}(Q) = \{(x^{-1}, x): x\in \Ac(P, Q), \, P\in \Pc\}.
\end{align}
Also, we denote by $x\mapsto \x = (x^{-1}, x)$ the canonical map
$\Ac \to \Ac^\tgt$.

If $T: \Ac\to \Bc$ is a morphism of quivers then
$\overline T: \Ac^{\tgt}\to \Bc^{\tgt}$ is the morphism
of loop bundles given by $\overline T(\x) = \overline{T(x)}$.

\bigbreak If $\Ac$, $\Bc$ are quivers over $\Pc$, then $\Ac
 \ftimes \Bc = \{(a,b) \in \Ac \times
\Bc: \tgt(a) = \sou(b)\}$ is a quiver over $\Pc$ with ${\mathfrak
s}(a,b) = \sou(a)$, $\tgt(a, b) = {\mathfrak e}(b)$. Thus, the
category $\Quiv(\Pc)$ of quivers over $\Pc$ is monoidal, with
$\otimes = \ftimes$ and  with unit object $(\Pc, \Pc, \id, \id)$.

\bigbreak This monoidal category does not seem to be symmetric in
any reasonable way. We have nevertheless two natural isomorphisms
playing the r\^ole of a ``weak" symmetry. A first one is the
natural isomorphism $\vartheta:  \Bc^{\op} {\,}_{\tgt}
\hspace{-0.1cm} \times_\sou \Ac^{\op} \to  \left(\Ac
{\,}_\tgt\hspace{-0.1cm} \times_\sou \Bc\right)^{\op}$ given by
\begin{equation}\label{vartheta}
\vartheta(y^{-1},x^{-1}) = (x,y), \qquad (x,y) \in \Ac
{\,}_\tgt\hspace{-0.1cm} \times_\sou \Bc.
\end{equation}
The second possibility is very similar. Let $\Bc {\,}_\sou
\hspace{-0.1cm} \times_\tgt \Ac = \{(b,a) \in \Bc \times \Ac:
\sou(b) = \tgt(a)\}$, a quiver over $\Pc$ with $\sou(b,a) =
\sou(a)$, $\tgt(b, a) = \tgt(b)$. Then we define $\tau: \Ac
{\,}_\tgt\hspace{-0.1cm} \times_\sou \Bc \to \Bc {\,}_\sou
\hspace{-0.1cm} \times_\tgt \Ac$ by
\begin{equation}\label{tau}
\tau(x,y) = (y,x), \qquad (x,y) \in \Ac {\,}_\tgt\hspace{-0.1cm} \times_\sou \Bc.
\end{equation}
These are related as follows. Let $\mu: \Bc {\,}_\sou
\hspace{-0.1cm} \times_\tgt \Ac  \to \left(\Bc^{\op} {\,}_{\tgt}
\hspace{-0.1cm} \times_\sou \Ac^{\op}\right)^{\op}$ be given by
$\mu(y,x) = (y^{-1}, x^{-1})$. Then the following diagram commutes:
\begin{equation*}
\xymatrix{ \Ac {\,}_\tgt\hspace{-0.1cm} \times_\sou \Bc
\ar[rr]^{\tau} & & \Bc {\,}_\sou \hspace{-0.1cm} \times_\tgt \Ac
\ar[dl]^{ \mu}
\\ & \ar[ul]^{\vartheta} \left(\Bc^{\op} \right. {\,}_{\tgt} \hspace{-0.1cm}
\times_\sou \left. \Ac^{\op}\right)^{\op}.  & }
\end{equation*}

\bigbreak
\subsection{Groupoids}

\

\bigbreak Let $\G$ be a groupoid  with base $\Pc$ and source and
end maps $\sou, \tgt: \G \rightrightarrows \Pc$. We identify $\Pc$
with a subset of $\G$ via the identity. We indicate the
composition $m(f, g)$  of two elements in a groupoid by
juxtaposition: $m(f, g) = fg$, and not $gf$. A group bundle is a
groupoid $\N$ with source = end; thus $\N = \coprod_{P\in \Pc}
\N(P)$.

\bigbreak A  morphism of groupoids $T: \G \to \K$ is a map
preserving the product; thus, it preserves also source and end,
and induces a map between the bases. If $\G$ and $\K$ have the
same base $\Pc$, we shall say that $T: \G \to \K$ is a morphism
of groupoids \emph{over} $\Pc$ if the restriction $\Pc \to \Pc$ is
the identity. Groupoids over $\Pc$ form a category
$\Gpd(\Pc)$.

\bigbreak A wide subgroupoid is a subgroupoid such that the
inclusion is a morphism of groupoids over $\Pc$. A subgroup bundle
of a groupoid is a subgroupoid that is a group bundle. There is a
largest wide subgroup bundle of a groupoid $\G$, namely
$\G^{\bundle} = \coprod_{P\in \Pc} \G(P)$; that is, we forget the
arrows between distinct points.

\bigbreak
If $(\N_i)_{i\in I}$ is a family of (wide) subgroupoids of a
groupoid $\G$ then $\bigcap_{i\in I} \N_i$, defined by
$\bigcap_{i\in I} \N_i(P; Q)$ $= \bigcap_{i\in I} \left(\N_i(P;
Q)\right)$, is a (wide) subgroupoid of $\G$.

\bigbreak A groupoid over $\Pc$ is a group object in the monoidal category $\Quiv (\Pc)$. Any groupoid being a quiver, we shall use all the
terminology above also for groupoids. Clearly, $\G(P)$ is a group
and it acts freely and transitively on the left on $\G(P, Q)$ and
on the right on $\G(Q, P)$, for any $P, Q\in \Pc$.

\bigbreak
Basic examples of groupoids are:
\begin{itemize}
\item A  group $G$, considered as the set of arrows of a category with a single object.
\item An equivalence relation $R$ on $\Pc$; $\sou$ and $\tgt$ are respectively the first
and the second projection, and the composition is given by
$(x,y)(y,v) = (x,v)$.
\end{itemize}

The equivalence relation where all the elements of $\Pc$ are
related is denoted $\Pc^2$ and called the \emph{coarse} groupoid
on $\Pc$.

\bigbreak
If $\G \rightrightarrows \Pc$ is any groupoid, then $\G \simeq
\coprod_{X \in \Pc/\conexo} \G_X$. Here $\G_X$ is the subgroupoid  on
the base $X$ defined by $\G_X(x,y) = \G(x,y)$, for all $x, y \in
X$. Furthermore, $\G_X = \G(x) \times X^2$ for any $x \in X$, $X \in \Pc/\conexo$. This description can be viewed as a structure theorem for groupoids.

\bigbreak Let $\N$ be a wide subgroup bundle of a groupoid $\G$
over $\Pc$. Then there are quivers $\N \backslash \G$, $\G / \N$,
equipped with surjective morphisms of quivers $\G \to \N
\backslash \G$ and $\G \to \G / \N$, defined by
$$(\N \backslash \G)(P, Q) = \N(P) \backslash \G(P, Q), \qquad
(\G / \N)(P, Q) = \G (P, Q)/ \N(Q), \qquad P, Q\in \Pc.$$

\bigbreak
We shall say that $\N$ is \emph{normal} if the following equivalent
conditions hold for any $P, Q\in \Pc$:
\begin{itemize}
\item For any $x\in \G(P, Q)$ and $n\in \N(Q)$, $xnx^{-1} \in \N(P)$.
\item For any $x\in \G(P, Q)$ and $n\in \N(Q)$, there exists $m \in \N(P)$
such that $xn = mx$.
\end{itemize}

\bigbreak  If $(\N_i)_{i\in I}$ is a family of normal subgroup
bundles of $\G$ then $\bigcap_{i\in I} \N_i$ is a normal subgroup
bundle.

\bigbreak If $\N$ is a normal wide subgroup bundle of $\G$ then
$\G / \N$ has groupoid multiplication, with the canonical map
$\pi: \G \to \G / \N$ being a morphism of groupoids.

\bigbreak Let $T: \G \to \K$ be a morphism of groupoids over a
$\Pc$. The \emph{kernel} of $T$ is the (wide and normal) subgroup
bundle
$$
\Ker T = \{g\in \G: T(g) \in \Pc\} = \coprod_{P\in \Pc} \Ker T(P), \qquad
\text{where } \Ker T(P) = \ker (T: \G(P) \to \K(P)).
$$
Let $f,g\in \G$. Then: $T(f) = T(g)$ iff there exists $n\in \Ker T$
with $f = ng$ iff there exists $m\in \Ker T$ with $f = gm$.
That is, $\Imm T \simeq \G / \Ker T$.

\bigbreak The largest subgroup bundle $\G^{\bundle}$ is clearly
normal;  the quotient $\G / \G^{\bundle}$ is the groupoid
associated to the equivalence relation $\approx$.

\bigbreak
\subsection{The free groupoid generated by a quiver}

\

\bigbreak It is natural to look for the left adjoint of the
obvious forgetful functor from $\Gpd(\Pc)$ to $\Quiv(\Pc)$. This
leads us to the construction of the free groupoid generated by a
quiver.

\begin{theorem} Let $\Ac$ be a quiver over $\Pc$. Then there
exists a groupoid $F(\Ac)$ over $\Pc$ provided with a morphism of
quivers $\iota: \Ac \to F(\Ac)$ satisfying the usual universal
property:

If $\G$ is a groupoid over $\Pc$ provided with a morphism of
quivers $\nu: \Ac \to \G$, then there is a unique map of groupoids
$\widehat\nu: F(\Ac) \to \G$ such that $\nu = \widehat \nu \,
\iota$.

The groupoid $F(\Ac)$ is unique up to isomorphisms with respect to this
property.
\end{theorem}

\pf  We first consider the quiver $\path (\D \Ac )$; its elements
will be called ``words in the alphabet $\Ac \cup \Ac^{\op}$". A
word $w = x_{1} x_{2} \dots x_{n}$, $x_i\in \D\Ac$,  is
\emph{reduced} if either $n = 0$, or there is no $i$ such that
$x_i = x_{i+1}^{-1}$.

Let $w = x_{1} x_{2} \dots x_{n}$, $x_i\in \D\Ac$,
be a word of length $n>0$ and
assume there exists $i$, $1\le i < n$, such that $x_i =
x_{i+1}^{-1}$. Then set $w' := x_{1} x_{2} \dots x_{i-1} x_{i + 2}
\dots x_{n}$, if $n >2$, or $w' := \id \sou(x_{1})$ if
$n=2$. The word $w'$ is called an \emph{elementary reduction} of $w$.
Furthermore, a word $\widetilde w$ is called a \emph{reduction} of $w$
if it can be attained from $w$ by a sequence of elementary reductions.

Reduction generates an equivalence relation in the usual way. Two
words $u$ and $v$ are \emph{equivalent}, denoted $u \sim v$, if
there is a sequence of words $w_{1}, w_{2} \dots, w_{N}$ with $N
\ge 1$, $u =w_{1}$, $v = w_{N}$ and either $w_i$ a reduction of
$w_{i+1}$, or $w_{i+1}$ a reduction of $w_{i}$, for all $i$, $1\le
i < N$. This is clearly an equivalence relation and the class of a
word $u$ is denoted $[u]$. Furthermore,

\begin{step}
In any class there is one and only one reduced word.
\end{step}

Let $w$ be a word. By a standard recurrence argument, there is at
least one reduced word in $[w]$, which is a reduction of $w$. To
prove the uniqueness, we consider the ``W-process" for a word. We
first set $\W_0 = \id \sou(w)$. If the length of $w$ is $n > 0$,
say $w = x_{1} x_{2} \dots x_{n}$, then we define recursively
\begin{align*}
\W_1 &= x_1,
\\ \W_{i+1} &= \begin{cases} &X, \qquad\quad \; \text{if }
\W_i \text{ is of the reduced form } Xx_{i+1}^{-1},
\\ &\W_ix_{i+1}, \quad \text{if } \W_i
\text{ is \emph{not} of the reduced form } Xx_{i+1}^{-1}.
\end{cases}
\end{align*}
Then $\W_0, \dots, \W_n$ are all reduced (by induction) and $\W_n
= w$ if $w$ is reduced. We have then a map $\path (\D \Ac ) \to
\{p\in \path (\D \Ac ): p$ is reduced$\}$, $w\mapsto \W_n$, which
is a retraction of the inclusion. We will now check that $w\sim u$
implies $\W_n = \Uc_m$, where $m$ is the length of $u$, and
$\Uc_0, \dots, \Uc_m$ is the W-process for $u$. In particular if
both $w$ and $u$ are reduced and equivalent, then necessarily
$w=u$.

\bigbreak So, assume that $w$ is an elementary reduction of $u =
x_{1} \dots x_{r} y y^{-1} x_{r+1} \dots x_{n}$, with $y\in \D
\Ac$. Clearly, $\Uc_0 = \W_0, \dots, \W_r = \Uc_r$. Now two cases
can happen:

a) $\W_r = \Uc_r$ is of the reduced form $Xy^{-1}$. Then $X$ is
\emph{not} of the reduced form $Yy$. Thus $\Uc_{r+1} = X$ and
$\Uc_{r+2} = Xy^{-1} = \W_r$. Hence $\Uc_{r+2 + i} = \W_{r + i}$,
$i\ge 0$.

b) $\W_r = \Uc_r$ is \emph{not} of the reduced form $Xy^{-1}$.
Then $\Uc_{r+1} = \Uc_r y$ and $\Uc_{r+2} = \Uc_r = \W_r$. Hence,
again, $\Uc_{r+2 + i} = \W_{r + i}$, $i\ge 0$.

This finishes the proof of the step.

\bigbreak The map $\path (\D \Ac ) {\,}_\tgt\hspace{-0.1cm}
\times_\sou \path (\D \Ac ) \to \path (\D \Ac )$,
$(x_1\dots x_n, y_1\dots y_m) \mapsto x_1\dots x_ny_1\dots y_m$,
if $n>0$, $m>0$; $(\id \sou(y_1), y_1\dots y_m) \mapsto y_1\dots y_m$,
etc., is called the \emph{juxtaposition}.
Let $F(\Ac) := \path (\D \Ac ) / \sim$.
\begin{step}
Juxtaposition induces a groupoid structure on $F(\Ac)$.
\end{step}

We first claim that juxtaposition descends to a map $\cdot:F(\Ac)
{\,}_\tgt\hspace{-0.1cm}\times_\sou F(\Ac) \to F(\Ac)$. We omit
the straightforward verification of the claim: ``$w\sim \tilde w$,
$u\sim \tilde u$ and $\tgt(w) = \sou(u)$ implies $wu\sim \tilde
w\tilde u$". Since $([u][v])[w]= [uvw] = [u]([v][w])$, $\cdot$ is
associative. The elements $[\id P]$, $P\in \Pc$, are partial
identities of the product $\cdot$. If $w = x_{1} x_{2} \dots
x_{n}$, then set $w^{-1} = x_{n}^{-1} x_{n-1}^{-1} \dots
x_{1}^{-1}$. Then $[w]^{-1} := [w^{-1}]$ is the inverse of $[w]$.
Thus $F(\Ac)$ is a groupoid.

\bigbreak
\begin{step}
The map $\iota: \Ac \to F(\Ac)$, $x \mapsto [x]$, is an injective
morphism of quivers and satisfies the required universal property.
\end{step}

Injectivity of $\iota$ follows from Step 1. If $\G$ is a groupoid
and $\nu: \Ac \to \G$ is a morphism of quivers, then $\nu$ can be
extended to $\Ac^{\op}$, then to $\path (\D \Ac )$ and finally to a
morphism of groupoids $\hat{\nu}:F(\Ac) \to \G$, which is easily
seen to be unique. \epf

\begin{obs}
This proof is an adaptation of the construction of the free group
generated by a set.
\end{obs}

\begin{exas}

\begin{itemize}
    \item The groupoid $F(\Ac)$ does not depend on the orientation
    of $\Ac$.
    \item If $\Pc$ has one element, then the groupoid $F(\Ac)$ is
    just the free group generated by $\Ac$.
    \item Suppose that $\Pc$ has exactly two elements $P$ and $Q$.
    If $\Ac$ consists of only one arrow from $P$ to $Q$, then
    $F(\Ac) \simeq \Pc^2$. If $\Ac$ consists only of arrows from $P$ to $Q$,
    say $\card \Ac(P,Q)=n$, then $F(\Ac) \simeq F^{n-1} \times \Pc^2$,
    where $F^{n-1}$ is the free group in $n-1$ variables.
\end{itemize}
\end{exas}

Let $\Ac$ be any quiver and let $\R$ be any subset of
$\coprod_{P\in \Pc} F(\Ac)(P)$. Then the \emph{groupoid presented
by $\Ac$ with relations $\R$} is the quotient of the  free
groupoid $F(\Ac)$ by the smallest normal wide subgroupoid
containing $\R$.

\bigbreak
Let $\G$ be a groupoid and $\Ac$, a sub-quiver of $\G$. The (wide)
\emph{subgroupoid generated by $\Ac$} is $<\Ac> :=$ the image
of the induced map $F(\Ac) \to \G$. In words, the elements of
$<\Ac>$ are compositions of elements of $\Ac$ or their inverses.
We say that
\emph{$\Ac$ generates $\G$} if  $<\Ac> = \G$.

\bigbreak
\subsection{Actions of groupoids}

\

\bigbreak Let $\G$ be a groupoid with base $\Pc$ and let $p: \Ec \to \Pc$ be
a map. A \emph{left action} of $\G$ on $p$ is a map $\fde : \G
{\,}_\tgt\hspace{-0.1cm}\times_p \Ec \to \Ec$ such that
\begin{equation}\label{axiomleft}
 p(g\fde e) = \sou(g), \qquad
 g \fde(h \fde e)  = gh \fde e, \qquad
\id \,{p(e)} \, \fde  e  =  e,
\end{equation}
for all $g, h \in \G$, $e \in \Ec$  composable appropriately.
Intertwiners of actions  of $\G$ on $p: \Ec \to \Pc$ and $p': \Ec'
\to \Pc$ are defined in the usual way.

Given a set $X$ there is an action of $\G$ on $p: \Pc\times X \to
\Pc$, where $p$ is the first projection, given by $g\fde (\tgt(g),
x) = (\sou(g), x)$ for all $x\in X$, $g\in \G$. An action on $p:
\Ec \to \Pc$ is \emph{trivial} if there exists a set $X$  and a
bijective intertwiner of actions $\nu: \Ec \to \Pc\times X$.

\begin{definition} Let $p: \Ec \to \Pc$ be a map. The groupoid $\aut p$,
or indistinctly $\aut \Ec$, is defined by
\begin{equation*}
\aut p =\left\{(P,x,Q): P,Q\in \Pc, \text{ and } x: \Ec_Q
\to \Ec_P \text{ is a bijection} \right\};
\end{equation*}
with source and end $\sou, \tgt : \aut p \to \Pc$ given by
$\sou(P,x,Q) = P$, $\tgt(P,x,Q) =Q$, $(P,x,Q) \in \aut p$; with
composition $(P,x,Q)(Q,y,R) = (P,xy, R)$, $(P,x,Q), (Q,y,R)\in
\aut p$; and with identities $\id P = (P,\id, P)$, $P\in \Pc$.
\end{definition}

Then there is an equivalence between left actions of $\G$ on $p$,
and morphisms of groupoids $\G \to \aut p$. Namely, if $\fde$ is a
left action, then the corresponding morphism $\rho: \G \to \aut p$
is $\rho(g) = (\sou (g), g \fde \underline{\;\;}\,, \tgt(g))$,
$g\in \G$.

\bigbreak Similarly, a \emph{right action} of $\G$ on $p$ is a map
$\fiz: \Ec {\,}_p \hspace{-0.1cm} \times_\sou \G \to \Ec$ such
that
\begin{equation}\label{axiomright}
p(e\fiz g) = \tgt(g), \qquad (e \fiz g) \fiz h  = e \fiz
gh, \qquad  e  \fiz \, \id \,{p(e)}  = e,
\end{equation}
for all $g, h \in \G$, $e \in \Ec$  composable appropriately. Left
and right actions are equivalent, by the rule $e\fiz g =
g^{-1}\fde e$.

\bigbreak
\subsection{Matched pairs of groupoids}\label{mpg}

\

\bigbreak It is convenient to introduce some alternative notation
for quivers and groupoids. We shall say that a groupoid
$\Vc$ is denoted \emph{vertically} if the source and end are named
respectively $t, b:\Vc \rightrightarrows \Pc$, where $t$ means
``top" and $b$ means ``bottom". The elements of $\Vc$ will  be
consequently depicted as vertical arrows going down. Similarly, a
groupoid $\Hc$ is denoted \emph{horizontally} if the source and
end are named respectively $l, r: \Hc \rightrightarrows \Pc$,
where $l$ means ``left" and $r$ means ``right"; the elements of
$\Hc$ will be depicted as horizontal arrows going right.

\begin{definition}\label{matchpairgpds}  \cite[Definition 2.14]{mk1}.
A \emph{matched pair of groupoids} is a pair of groupoids $(\Vc,
\Hc)$ over $\Pc$ with $\Vc$ denoted vertically and $\Hc$
horizontally, endowed with a left action $\fde : \Hc {\,}_r
\hspace{-0.1cm} \times_t \Vc \to \Vc$ of $\Hc$ on $t: \Vc \to
\Pc$, and a right action $\fiz : \Hc {\,}_r \hspace{-0.1cm}
\times_t \Vc \to \Hc$ of $\Vc$ on $r: \Hc \to \Pc$, satisfying
\begin{flalign}
\label{mp-0.7}  &b(x\fde g) = l(x \fiz g),&
\\ \label{mp-3} & x \fde fg  = (x \fde f) ((x \fiz f) \fde g), &
\\ \label{mp-4} & xy \fiz g = (x \fiz (y \fde g)) (y \fiz g), &
\end{flalign}
for all $f, g \in \Vc$, $x, y \in \Hc$ such that the compositions
are possible.
\end{definition}

Given $P\in \Pc$, there are two identities: $\id_{\Hc} P \in \Hc$
and $\id_{\Vc} P\in \Vc$. We omit the subscript unless some
emphasis is needed.

\begin{definition}\label{exfactgpds} \cite{mk1}.
Let $\D \rightrightarrows \Pc$ be a groupoid. An  \emph{exact
factorization} of $\D$ is a pair of wide subgroupoids $\Vc$,
$\Hc$, such that the multiplication map $\Vc {\,}_b\times_l \Hc
\to \D$ is a bijection.
\end{definition}

We shall not give the complete definition of a vacant double
groupoid (due to Ehresmann),  see \cite{AN} for a full discussion and
some historical references. Informally, a
double groupoid is a set of boxes with two partial compositions, a
horizontal one and a vertical one. Each box has horizontal and
vertical sides; the set of all horizontal sides form a groupoid on
their own, and the same holds for the vertical sides. These
two groupoids share the set of points, which are the corners of
the boxes. In short, a double groupoid is a collection of sets and maps
$$\begin{matrix} \qquad\Bg &\overset{t,b}\rightrightarrows &\Hc \qquad
\\ l,r \downdownarrows &\qquad&\downdownarrows l,r
\\ \qquad \Vc &\underset{t,b}\rightrightarrows &\Pc \qquad\end{matrix}$$
such that all four sides in this diagram are groupoids, and satisfying
some compatibility conditions.

\bigbreak
A double groupoid is vacant if any pair of a horizontal
and a vertical side with a common point determines a unique box.

\bigbreak We recall some notation needed later. If $X$ is a box
then $X^{-1}$ is the box obtained inverting the horizontal
and the
vertical arrows. If $g\in \Vc$ with $t(g) = P$ and $b(g) = Q$ then
the horizontal identity of $g$ is the box $$\idd \, g =
\begin{matrix} \quad \id_{\Hc} P  \quad \\  g \, \boxee \,\, g
\\ \quad  \id_{\Hc} Q   \quad
\end{matrix}.$$

The vertical identity of a horizontal arrow is defined similarly.

\begin{proposition}\label{equiv-matchedpair} \cite[Theorems 2.10 and 2.15]{mk1}
The following notions are equivalent.

\begin{enumerate}
\item Matched pairs of groupoids.
\item Groupoids with an exact factorization.
\item Vacant double groupoids.
\end{enumerate}
\end{proposition}

We sketch the main parts of the correspondence needed later. See
\cite[Theorems 2.10 and 2.15]{mk1} and also  \cite[Prop 2.9]{AN}
for more details.

\pf If ($\Vc$, $\Hc$) is a matched pair of groupoids, then the
\emph{diagonal} groupoid $\Vc\bowtie \Hc\rightrightarrows \Pc$ is
defined as follows: $\Vc\bowtie \Hc := \Vc {\,}_b \hspace{-0.1cm}
\times_l \Hc$, with composition $(f,y) (h, z) = (f (y\fde h),
(y\fiz h) z)$, with source $\sou: \Vc {\,} \to \Pc$, $\sou(f, y) =
t(f)$, and with end $\tgt:  \Vc {\,} \to \Pc$, $\tgt(f, y) =
r(y)$. Clearly, $\Vc$ and $\Hc$ can be identified with
subgroupoids of $\Vc \bowtie \Hc$ forming an exact factorization.

\bigbreak Conversely, let $\Vc$, $\Hc$ be an exact factorization
of a groupoid $\D$; that is, for any $\alpha \in \D$, there exist
unique $f\in \Vc$, $y \in \Hc$, such that $\alpha = fy$. Let
$\begin{CD} \Hc @<\fiz<<  \Hc {\,}_r \hspace{-0.1cm} \times_t \Vc
@>\fde>>  \Vc
\end{CD}$ be given by $xg = (x\fde g)(x\fiz g)$, $(x,g) \in \Hc
{\,}_r \hspace{-0.1cm} \times_t \Vc$. Then $\Vc$, $\Hc$, together
with these actions, form a matched pair.

\bigbreak Similarly, if ($\Vc$, $\Hc$) is a matched pair of
groupoids  then we define $\Bg := \Hc {\,}_r\times_t \Vc$. We
represent $X = (x,g) \in \Hc {\,}_r\times_t \Vc$ by $$X =
\, \begin{matrix} \qquad\quad  x  \quad \\  x \fde g \, \boxee \,\, g
\\ \qquad\quad  x \fiz g  \quad
\end{matrix}.$$ Then $\T = (\Bg, \Vc, \Hc, \Pc)$ is a vacant double groupoid.
For later use, we record the description of the groupoid $\Bg
\rightrightarrows \Vc$: $l(x,g) = x\fde g$, $r(x,g) = g$,
$(x,g)(y,h) = (xy, h)$ if $g = y\fde h$.
 The construction is reversible and gives the opposite
implication. \epf

\bigbreak
\subsection{Semidirect products}

\

\bigbreak Let $\Vc$ be a groupoid denoted vertically. Assume that
$\Hc = \N$ is a group bundle, with $p := l = r$. The
\emph{trivial} action of $\N$ on $\Vc$ is $x\fde g = g$, $(x,g)
\in \N  {\,}_p \hspace{-0.1cm} \times_t \Vc.$ Also, a right action
of $\Vc$ on $\N$ is \emph{by group bundle automorphisms} if $
xy\fiz g = (x\fiz g)(y\fiz g)$, $x,y \in \N$, $g \in  \Vc, $
composable.

\bigbreak It is easy to see that a right action is a compatible
with the trivial action iff it is by group bundle automorphisms.
We shall  denote the corresponding diagonal groupoid by $\Vc
\ltimes \N$ and call it a \emph{semidirect product}. The
projection $\Vc \ltimes \N \to \Vc$ is a morphism of groupoids,
and has a section of groupoids $S: \Vc \to \Vc \ltimes \N$, $S(f)
= (f, \id b(f))$.

\bigbreak
Conversely, let $T: \G \to \K$ be a morphism of groupoids over $\Pc$.
Then $\G$ acts on the kernel $\N$ of $T$ by the adjoint action:
$n\fiz g = g^{-1}ng$. If there is a section $S: \K \to \G$, then
$\G \simeq \K \ltimes \N$.

\bigbreak The structure theorem of groupoids can be phrased in
this language  as the isomorphism of groupoids $\G \simeq \R
\times \G^{\bundle}$ where $\R$ is the groupoid associated to the
equivalence relation $\approx$.

\bigbreak
\subsection{Actions of matched pairs of groupoids}

\

\bigbreak Let ($\Vc$, $\Hc$) be a matched pair of groupoids over $\Pc$.

\begin{definition}\label{tquiver} \cite{AA}.
A (set-theoretic) {\em representation} of
$(\Vc,\Hc)$ is a triple $(\Ac, \fde, \vert \ \vert)$, where
\begin{itemize}
\item $\Ac$ is a quiver over $\Pc$,
\item $\fde : \Hc {\,}_r \times_\sou\Ac \to \Ac$ is a left action
of $\Hc$ on $\sou$, and
\item $\vert \ \vert: \Ac \to \Vc$ is a morphism of quivers over
$\Pc$, called the grading, such that
\end{itemize}
\begin{equation}\label{compcond}
\vert x \fde a \vert = x \fde \vert a \vert, \qquad (x, a) \in \Hc {\,}_r \times_\sou \Ac.
\end{equation}
We shall say simply ``$\Ac$ is a representation of $(\Vc,\Hc)$".

\bigbreak Morphisms of representations of $(\Vc,\Hc)$ are
morphisms of quivers intertwining the actions of $\Hc$ and
preserving the grading $\vert\,\vert$. Thus representations of
$(\Vc,\Hc)$ form a category $\Rep(\Vc,\Hc)$; this is a monoidal
subcategory of $\Quiv(\Pc)$. Namely, if  $\Ac$ and $\Bc$ are two
representations of $(\Vc,\Hc)$, then $\Ac \ftimes \Bc$ is also a
representation of $(\Vc,\Hc)$, with respect to the action and
grading given by
\begin{align}\label{action-pt}
x \fde (a, b)  &= (x \fde a,  (x \fiz  \vert a \vert) \fde b),
\\\label{grado-pt} \vert (a,b)\vert &= \vert a\vert  \vert b\vert,
\end{align}
$x\in \Hc$, $(a,b) \in \Ac \ftimes \Bc$. See \cite{AA} for details.
\end{definition}

\begin{obs} We can state the preceding notion in terms of the
associated vacant double groupoid. A representation of $(\Vc,\Hc)$
is the same as a left action of the groupoid $\Bg
\rightrightarrows \Vc$. In fact, we have:

(a). Let $\Vc$ be a quiver over $\Pc$ and let $p: \Ec \to \Vc$ be
a map. Then $\Ec$ is a quiver over $\Pc$ with source $\sou\circ p$
and end $\tgt\circ p$. Moreover $p$ is a morphism of quivers. We
have a functor from the category of sets over $\Vc$ to
$\Quiv(\Pc)$.

(b). Let $\Vc$ be a groupoid over $\Pc$. Then the category of sets
over $\Vc$ is monoidal; if $p: \Ec \to \Vc$ and $q: \F \to \Vc$
are maps then the tensor product is $\Ec \ftimes \F$, with grading
\eqref{grado-pt}. The functor in (a) is monoidal too.

(c). Let ($\Vc$, $\Hc$) be a matched pair of groupoids over $\Pc$.
If $\Ac$ is a representation of $(\Vc,\Hc)$ then we define a left
action of the groupoid $\Bg \rightrightarrows \Vc$ on the map
$\vert \ \vert: \Ac \to \Vc$ by the rule
$$
(x,g) \fde a = x\fde a, \qquad \text{ if } g = \vert a \vert.
$$
Conversely, a left action of the groupoid $\Bg \rightrightarrows
\Vc$ on a map $\vert \ \vert: \Ac \to \Vc$ defines a
representation of $(\Vc,\Hc)$ by the same rule. We omit the
straightforward verifications.

However there is no obvious translation of the monoidal structure
on $\Rep(\Vc, \Hc)$ to the actions of $\Bg \rightrightarrows \Vc$.
In other words we can not extend the definition of representations
of matched pairs to arbitrary double groupoids.
\end{obs}

\begin{lema}\label{repndual} If $\Ac$ is a representation of
$(\Vc,\Hc)$ then $\Ac^{\op}$ is a representation of $(\Vc,\Hc)$ with
\begin{equation} x \fde a^{-1} = \big((x \fiz \vert a\vert^{-1}) \fde a\big)^{-1}, \qquad \vert a^{-1}\vert = \vert a\vert^{-1}, \qquad
x\in \Hc, a\in \Ac.
\end{equation}
\end{lema}

\pf Straightforward, using formula (1.17) in \cite{AA}. \epf

Since disjoint unions of representations  are again
representations, we conclude that $\D \Ac$ is a representation of $(\Vc,\Hc)$.

\section{The quiver-theoretical  braid equation}

In this section we introduce our main object of study-- the quiver
theoretical braid equation. We establish some basic properties and
present methods to construct examples.

\bigbreak
\subsection{The quiver-theoretical  quantum Yang-Baxter equation}

\

\bigbreak Let $\Ac$ be a quiver and let
$\sigma: \Ac{\,}_\tgt\hspace{-0.1cm} \times_\sou \Ac
\to \Ac{\,}_\tgt\hspace{-0.1cm} \times_\sou
\Ac$ be an isomorphism of quivers. We set $\Ac^n := \Ac
{\,}_\tgt\hspace{-0.1cm} \times_\sou \Ac {\,}_\tgt\hspace{-0.1cm}
\times_\sou \dots {\,}_\tgt\hspace{-0.1cm} \times_\sou \Ac
{\,}_\tgt\hspace{-0.1cm} \times_\sou \Ac$, $n$-times; and
$\sigma_{i, i+1} := \id_{\Ac^{i-1}} \times\sigma \times
\id_{\Ac^{n-i-1}}: \Ac^n \to \Ac^n$, an isomorphism of quivers.

\bigbreak A \emph{solution of the quiver-theoretical braid equation over
$\Pc$} (a \emph{``solution"} or a \emph{``braided quiver"}, for short)
is a pair formed by a quiver
$\Ac$ and an isomorphism of quivers $\sigma:  \Ac \ftimes \Ac \to  \Ac \ftimes \Ac$
such that
\begin{equation}
\label{qtBE}  (\sigma \times \id) (\id \times \sigma)(\sigma
\times \id) = (\id \times \sigma)(\sigma \times \id)(\id \times
\sigma):
\Ac {\,}_\tgt \hspace{-0.1cm}\times_\sou \Ac {\,}_\tgt \hspace{-0.1cm}\times_\sou  \Ac
\to \Ac {\,}_\tgt \hspace{-0.1cm}\times_\sou \Ac {\,}_\tgt
\hspace{-0.1cm}\times_\sou  \Ac.
\end{equation}

A solution $(\Ac, \sigma)$ is a \emph{symmetry} when $\sigma^2 =
\id$. If $(\Ac, \sigma)$ is a solution, then the braid group
$\mathbb B_n$ acts by automorphisms of quivers  on $ \Ac^n$ for
any $n\ge 2$; if $(\Ac, \sigma)$ is a symmetry, this action
descends to an action of the symmetric group $\mathbb S_n$.

\bigbreak Analogously, a \emph{solution of the quiver-theoretical
quantum Yang-Baxter equation over $\Pc$} is a pair formed by a
quiver $\Ac$  and an isomorphism of quivers $R:  \Ac {\,}_\tgt
\hspace{-0.1cm}\times_\sou \Ac \to  \Ac {\,}_\sou \hspace{-0.1cm}
\times_\tgt \Ac$ such that
\begin{equation}
\label{qtQYBE} R_{12} R_{13} R_{23} = R_{23}R_{13}R_{12}:
 \Ac {\,}_\tgt \hspace{-0.1cm}\times_\sou \Ac {\,}_\tgt
 \hspace{-0.1cm}\times_\sou \Ac \to \Ac {\,}_\sou \hspace{-0.1cm}\times_\tgt
 \Ac {\,}_\sou \hspace{-0.1cm}\times_\tgt  \Ac.
\end{equation}
Note that $R$ and both members of the equality \eqref{qtQYBE} are
isomorphisms between \emph{different} quivers. In fact, the
precise  meaning of the members of \eqref{qtQYBE} is
given by the commutativity of the diagrams:
$$
\begin{CD}  \Ac {\,}_\tgt \hspace{-0.1cm}\times_\sou \Ac {\,}_\tgt
\hspace{-0.1cm}\times_\sou \Ac @>R_{23}R_{13}R_{12}>> \Ac {\,}_\sou
\hspace{-0.1cm}\times_\tgt \Ac {\,}_\sou \hspace{-0.1cm}\times_\tgt  \Ac \\
@V R_{12}VV @AA R_{23}A \\
(\Ac {\,}_\sou \hspace{-0.1cm}\times_\tgt \Ac) {\,}_\tgt \hspace{-0.1cm}
\times_\sou \Ac @>R_{13}>> \Ac {\,}_\sou \hspace{-0.1cm}\times_\tgt
(\Ac {\,}_\tgt \hspace{-0.1cm}\times_\sou  \Ac),
\end{CD} \qquad
\begin{CD}  \Ac {\,}_\tgt \hspace{-0.1cm}\times_\sou \Ac {\,}_\tgt
\hspace{-0.1cm}\times_\sou \Ac @>R_{12} R_{13} R_{23}>> \Ac {\,}_\sou
\hspace{-0.1cm}\times_\tgt \Ac {\,}_\sou \hspace{-0.1cm}\times_\tgt  \Ac \\
@V R_{23}VV @AA R_{12}A \\
\Ac {\,}_\tgt \hspace{-0.1cm}\times_\sou (\Ac {\,}_\sou \hspace{-0.1cm}
\times_\tgt \Ac) @>R_{13}>> (\Ac {\,}_\tgt \hspace{-0.1cm}\times_\sou \Ac)
{\,}_\sou \hspace{-0.1cm}\times_\tgt  \Ac.
\end{CD}
$$

\bigbreak As usual, there is a bijective correspondence between
solutions of the quiver-theoretical quantum Yang-Baxter equation
and solutions of the quiver-theoretical braid equation over $\Pc$.
Namely, let $\sigma:  \Ac {\,}_{\mathfrak e} \hspace{-0.1cm}
\times_\sou \Ac \to  \Ac {\,}_\tgt \times_\sou \Ac$ be an
isomorphism of quivers, let $\tau: \Ac {\,}_\tgt
\hspace{-0.1cm}\times_\sou \Ac \to  \Ac {\,}_\sou \hspace{-0.1cm}
\times_\tgt \Ac$ be given by \eqref{tau} and let $R:= \tau\sigma$.
Then $R$ is a solution of the quiver-theoretical quantum
Yang-Baxter equation if and only if $\sigma$ is a solution of the
quiver-theoretical braid equation.

\bigbreak
\subsection{Non-degenerate  solutions}

\

\bigbreak Let $(\Ac, \sigma)$ be a braided quiver.  We define maps
$\fde, \fiz, \fd, \fz: \Ac {\,}_\tgt \hspace{-0.1cm}\times_\sou
\Ac \to \Ac$ by
\begin{align}\label{trenzas-notac}
\sigma(x, y) &=  (x\fde y, x \fiz y),
\\ \label{trenzas-notac1}
\sigma^{-1}(x, y) &=  (x\fd y, x \fz y),
\end{align}
$(x, y) \in \Ac \ftimes \Ac$. Clearly,
\begin{equation}\label{esquinas-braided}
\mathfrak s(x) =  \mathfrak s(x\fde y), \qquad \mathfrak e(x\fde
y) = \mathfrak s(x \fiz y), \qquad \mathfrak e(x \fiz y) =
\mathfrak e(y), \qquad (x, y) \in \Ac {\,}_\tgt \hspace{-0.1cm}
\times_\sou \Ac.
\end{equation}
We shall use below the relations between $\fde, \fiz, \fd, \fz$, like
\begin{equation}\label{trenzas-notac1.5}
x = (x\fde g) \fd (x\fiz g), \quad g = (x\fde g) \fz (x\fiz g), \quad y
= (y\fd h) \fde (y\fz h), \quad h = (y\fd h) \fiz (y\fz h),
\end{equation}
for composable $x,g,y,h$.

\bigbreak
The braid equation \eqref{qtBE} can be restated as
$(\id \times \sigma)(\sigma \times \id) (\id \times\sigma)^{-1}=
(\sigma \times \id)^{-1}(\id \times \sigma)(\sigma \times \id)$.
The equality of the first two components in this identity, specialized at $(h, f, u)\in  \Ac\ftimes \Ac\ftimes \Ac$, gives
\begin{align}\label{tecnica1}
h \fde \left( f\fd u \right) &=
(h\fde f) \fd \left[(h\fiz f)\fde u \right],
\\ \label{tecnica2}
\left[h \fiz (f \fd u) \right] \fde \left( f\fz u \right) &=
(h\fde f) \fz \left[(h\fiz f)\fde u \right].
\end{align}

\bigbreak
\begin{obs}\label{rigid} Let $\C$ be any tensor category.
Then we can define a solution of the braid equation in $\C$ as a
pair $(V, c)$, where $V$ is an object of $\C$ and $c:V \otimes V
\to V \otimes V$ is an invertible arrow satisfying $(c\otimes
\id)(\id\otimes c)(c\otimes \id) = (\id\otimes c)(c\otimes
\id)(\id\otimes c)$. Assume further that $\C$ is \emph{rigid}.
Then we say that a solution $(V, c)$ is \emph{rigid} if the map
$c^{\flat}:V^*\otimes V\to V\otimes V^*$ is invertible, where
$$c^{\flat}=(\ev_V\otimes\id_{V\otimes V^*}) (\id_{V^*}\otimes
c\otimes\id_{V^*}) (\id_{V^*\otimes V}\otimes \coev_V).$$ The
category $\Quiv(\Pc)$ is not rigid but the analog
of rigid solutions is given by the following definition.
\end{obs}

\begin{definition}
A solution $(\Ac, \sigma)$ is \emph{non-degenerate} if
\begin{equation}\label{trenzas-notac2}
x \fde \underline{\;\;}: \sou^{-1}(\tgt (x)) \to
\sou^{-1}(\sou (x)) \text{ and } \underline{\;\;} \fiz x: \tgt^{-1}(\sou (x))
\to \tgt^{-1}(\tgt (x))
\end{equation}
are bijections for any $x\in \Ac$.

\bigbreak
We next introduce the structure groupoid of a solution, a generalization of
definitions in \cite{ess,lyz1,s}. It plays the r\^ole of the FRT-bialgebra in this context.

\begin{definition}\label{frt} The \emph{structure groupoid} of the braided
quiver $(\Ac, \sigma)$ is the groupoid $\gstr$ generated by $\Ac$ with
relations
\begin{equation}
xy =  (x\fde y)(x \fiz y), \qquad (x, y) \in \Ac {\,}_{\mathfrak
e} \hspace{-0.1cm}\times_\sou \Ac.
\end{equation}\label{str-gpd-def1}
Equivalently, it can be defined as the groupoid
generated by $\Ac$ with relations
\begin{equation}\label{str-gpd-def2}
xy =  (x\fd y)(x \fz y), \qquad (x, y) \in \Ac {\,}_{\mathfrak
e} \hspace{-0.1cm}\times_\sou \Ac.
\end{equation}
\end{definition}

\bigbreak Note that if $\fde: \Ac {\,}_\tgt
\hspace{-0.1cm}\times_\sou \Ac \to \Ac$ is a map such that $x \fde
\underline{\;\;}: \sou^{-1}(\tgt (x)) \to \sou^{-1}(\sou (x))$ is
a bijection for any $x\in \Ac$ then it extends to a left action of
the free groupoid $F(\Ac)$ on $\sou: \Ac \to \Pc$. Similarly, a
map $\fiz$ with the analogous property induces a right action of
$F(\Ac)$ on $\tgt: \Ac \to \Pc$.

\begin{lema}\label{accion1} Let $\Ac$ be a quiver and $\sigma:  \Ac {\,}_\tgt
\hspace{-0.1cm}\times_\sou \Ac \to  \Ac {\,}_\tgt \hspace{-0.1cm}
\times_\sou \Ac$ be an isomorphism of quivers. Then $(\Ac,
\sigma)$ is a non-degenerate solution if and only if
\begin{flalign}
\label{cond-structur1} & \fde \text{ extends to a left action of }
\gstr \text{ on } \sou: \Ac \to \Pc;& \\
\label{cond-structur2} & \fiz \text{ extends to a right action of
} \gstr \text{ on } \tgt: \Ac \to \Pc;& \\
\label{cond-structur3} & (x\fde y) \fiz ((x\fiz y)\fde z) = (x
\fiz (y\fde z)) \fde (y\fiz z), \quad \text{for all } (x,y,z) \in
\Ac {\,}_\tgt \hspace{-0.1cm}\times_\sou \Ac {\,}_\tgt
\hspace{-0.1cm}\times_\sou \Ac.&
\end{flalign}
\end{lema}
\pf We fix $(x,y,z) \in  \Ac {\,}_\tgt \hspace{-0.1cm}
\times_{\mathfrak s} \Ac \ftimes
\Ac$ and compute:
\begin{align*}
(\sigma \times \id) (\id \times \sigma)(\sigma \times \id) (x,y,z)
&=
((x \fde y) \fde \left( (x \fiz y) \fde z\right), (x\fde y) \fiz ((x\fiz y)\fde z),
 (x\fiz y) \fiz z ),\\
(\id \times \sigma)(\sigma \times \id)(\id \times \sigma) (x,y,z)
&= (x \fde(y \fde z),  (x \fiz (y\fde z)) \fde (y\fiz z), (x \fiz
(y\fde z)) \fiz (y\fiz z)).
\end{align*}
Then $(\Ac, \sigma)$ is a solution iff 3 equalities hold, the
second one being \eqref{cond-structur3}. In presence of
non-degeneracy, the first of these equalities is equivalent to
\eqref{cond-structur1} and the third to \eqref{cond-structur2}.
\epf

If $(\Ac, \sigma)$ is a non-degenerate solution then
we shall denote by
\begin{equation}\label{trenzas-notac3}
h^{-1} \fde \underline{\;\;}, \text{ respectively } \underline{\;\;}\,\fiz
g^{-1}, \text{ the inverse of } h \fde \underline{\;\;}\,,
\text{ respectively } \underline{\;\;}\,\fiz g.
\end{equation}
\end{definition}

\bigbreak Let $(\Ac, \sigma)$ be a solution. Then $(\Ac,
\sigma^{-1})$ is a solution. Furthermore, if $\vartheta:
\Ac^{\op} {\,}_{\tgt} \hspace{-0.1cm} \times_\sou \Ac^{\op} \to
\Ac {\,}_\tgt\hspace{-0.1cm} \times_\sou \Ac$ is given by
\eqref{vartheta},  then $(\Ac^{\op},
\vartheta^{-1}\sigma\vartheta)$, as well as $(\Ac^{\op},
\vartheta^{-1}\sigma^{-1}\vartheta)$, are solutions. We denote
$\sigma^* = \vartheta^{-1}\sigma^{-1}\vartheta$ and
$\sigma^*(x^{-1}, y^{-1}) =  (x^{-1}\fde y^{-1}, x^{-1} \fiz
y^{-1})$; thus
\begin{equation}\label{trenzas-notac4}
x^{-1}\fde y^{-1} = (y \fz x)^{-1}, \qquad
x^{-1} \fiz y^{-1} = (y \fd x)^{-1}.
\end{equation}

\bigbreak
Let $\gstro$ be the structure groupoid of $(\Ac^{\op}, \sigma^{*})$.
Then the map $\Ac \to \Ac^{\op}$, $x\mapsto x^{-1}$, induces an isomorphism
of groupoids
$\gstr \to \gstro$, \emph{cf.} \eqref{str-gpd-def1}.
We shall identify $\gstr = \gstro$ via this isomorphism in what follows.

\begin{lema}\label{rigid-dual} Let $(\Ac, \sigma)$ be a non-degenerate
solution. Then

(a) $(\Ac^{\op}, \sigma^{*})$ is non-degenerate.

(b) $(\Ac, \sigma^{-1})$ is non-degenerate.

(c).  We have in $\gstr$:
\begin{align} \label{gstr-aux1}
xy^{-1} &= \left(x \fde y^{-1} \right)\left(x \fiz y^{-1} \right),
\\ \label{gstr-aux2}
x^{-1}z &= \left(x^{-1} \fde z\right)\left(x^{-1} \fiz z\right),
\end{align}
$x,y,z\in \Ac$, $\sou (x) = \sou(z)$, $\tgt(x) = \tgt(y)$.
\end{lema}

By (a), there are actions $\fde: \gstr {\,}_\tgt \hspace{-0.1cm}\times_
\sou \Ac^{\op} \to \Ac^{\op}$, $\fiz: \Ac^{\op} {\,}_\tgt \hspace{-
0.1cm} \times_\sou \gstr \to \Ac^{\op}$.

\pf (a). We define $\fde: \Ac \ftimes
\Ac^{\op} \to \Ac^{\op}$, $\fiz: \Ac^{\op} {\,}_\tgt
\hspace{-0.1cm}\times_\sou \Ac \to \Ac^{\op}$ by
\begin{equation}\label{trenzas-notac5}
x \fde g^{-1} = ((x\fiz g^{-1}) \fde g)^{-1}, \qquad h^{-1} \fiz y
= (h\fiz (h^{-1} \fde y))^{-1}.
\end{equation}
We compute
\begin{align*}
x^{-1}\fde (x \fde g^{-1}) &= x^{-1}\fde ((x\fiz g^{-1}) \fde g)^{-1}
= \left( ((x\fiz g^{-1}) \fde g) \fz x \right)^{-1}
\\ &= \left( \left((x\fiz g^{-1}) \fde g \right) \fz
\left(\left(x \fiz g^{-1}\right) \fiz g \right) \right)^{-1}
= g^{-1},\end{align*}
and also
$$x\fde (x^{-1} \fde g^{-1}) = x\fde (g\fz x)^{-1}
= \left((x\fiz (g\fz x)^{-1}) \fde (g\fz x) \right)^{-1}
= \left((x\fd g) \fde (g\fz x) \right)^{-1} =g^{-1}.$$
This means that
$x^{-1} \fde \underline{\;\;}\,: \Ac^{\op} \to \Ac^{\op}$ is the partial
inverse of $x \fde \underline{\;\;}\,: \Ac^{\op} \to \Ac^{\op}$. Similarly
for the right action, and (a) is proved. Then (b) follows from (a) by
\eqref{trenzas-notac4}.

(c).  By definition of $\gstr$, $(x\fiz y^{-1}) y
= \left((x\fiz y^{-1}) \fde y\right)\left((x\fiz y^{-1}) \fiz y\right)
= \left((x\fiz y^{-1}) \fde y\right)x$, thus
$xy^{-1}
= \left((x\fiz y^{-1}) \fde y\right)^{-1} \left(x \fiz y^{-1} \right)
= \left(x \fde y^{-1} \right)\left(x \fiz y^{-1} \right),$
proving \eqref{gstr-aux1}. The proof of \eqref{gstr-aux2} is similar.
\epf

\bigbreak
Putting together \eqref{trenzas-notac2}, \eqref{trenzas-notac3},
\eqref{trenzas-notac4} and \eqref{trenzas-notac5}, we have maps
$\fde, \fiz: \D \Ac {\,}_\tgt \hspace{-0.1cm}\times_\sou \D \Ac \to
\D \Ac$. We can then define
$\osigma: \D \Ac {\,}_\tgt \hspace{-0.1cm}\times_\sou \D \Ac \to
\D \Ac {\,}_\tgt \hspace{-0.1cm}\times_\sou \D \Ac$ by
\begin{equation}\label{trenzas-notac6}
\osigma(x, y) =  (x\fde y, x \fiz y),\qquad (x, y) \in \D \Ac \ftimes \D\Ac.
\end{equation}

\begin{lema}\label{rigid-double}  If $(\Ac, \sigma)$ is a
non-degenerate solution then $(\D \Ac, \osigma)$ is a
non-degenerate solution.
\end{lema}

\pf  Once the validity of the braid equation is established, the
rigidity will be clear by construction. Now it is necessary to
verify the equality \eqref{qtBE} on 8 subsets
of $\D\Ac \ftimes \D\Ac$. It holds in
$\Ac {\,}_\tgt \hspace{-0.1cm}\times_\sou \Ac {\,}_\tgt \hspace{-0.1cm}\times_\sou  \Ac$
by hypothesis and it holds in
$\Ac^{\op} {\,}_\tgt \hspace{-0.1cm}\times_\sou \Ac^{\op} {\,}_\tgt \hspace{-0.1cm}
\times_\sou  \Ac^{\op}$ by Lemma \ref{rigid-dual} (a).
Since $\left((\Ac^{\op})^{\op}, (\sigma^{*})^{*}\right)
= (\Ac, \sigma)$ we are reduced to three cases.

\bigbreak
\emph{Case I.} $(\sigma \times \id) (\id \times \sigma)(\sigma
\times \id) \overset{?}= (\id \times \sigma)(\sigma \times \id)(\id \times \sigma):
\Ac {\,}_\tgt \hspace{-0.1cm}\times_\sou \Ac^{\op} {\,}_\tgt
\hspace{-0.1cm}\times_\sou  \Ac^{\op}
\to \Ac^{\op} {\,}_\tgt \hspace{-0.1cm}\times_\sou \Ac^{\op} {\,}_\tgt
\hspace{-0.1cm}\times_\sou  \Ac.$

\bigbreak
We fix
$(x,y^{-1},z^{-1}) \in  \Ac {\,}_\tgt \hspace{-0.1cm}\times_\sou \Ac^{\op} {\,}_\tgt \hspace{-0.1cm}\times_\sou  \Ac^{\op}
$ and compute both sides of the desired equality:

\begin{align*}
LHS &=
((x \fde y^{-1}) \fde \left( (x \fiz y^{-1}) \fde z^{-1}\right),
(x\fde y^{-1}) \fiz ((x\fiz y^{-1})\fde z^{-1}),
 (x\fiz y^{-1}) \fiz z^{-1} ),\\
RHS &= (x \fde(y^{-1} \fde z^{-1}),
(x \fiz (y^{-1}\fde z^{-1})) \fde (y^{-1}\fiz z^{-1}),
(x \fiz(y^{-1}\fde z^{-1})) \fiz (y^{-1}\fiz z^{-1})).
\end{align*}

The first components of the LHS and the RHS are equal by \eqref{gstr-aux1}, while the last are equal by definition of $\gstro$.
It remains to show the equality of the middle components:

\begin{equation}\label{rigid-double-aux1}
(x\fde y^{-1}) \fiz ((x\fiz y^{-1})\fde z^{-1})
\overset{?}= (x \fiz (y^{-1}\fde z^{-1})) \fde (y^{-1}\fiz z^{-1}).
\end{equation}
Now RHS of \eqref{rigid-double-aux1} $= \left[x^{-1} \fiz \left( x \fde (y^{-1}\fde z^{-1})\right) \right]^{-1} \fde (y^{-1}\fiz z^{-1})$ by \eqref{trenzas-notac5}. Thus, we are reduced to prove

\begin{align*} y^{-1}\fiz z^{-1} &\overset{?}=
\left[x^{-1} \fiz \left( x \fde (y^{-1}\fde z^{-1})\right) \right] \fde \left((x\fde y^{-1}) \fiz ((x\fiz y^{-1})\fde z^{-1}) \right)
\\ &= \left[x^{-1} \fiz \left( (x\fde y^{-1}) \fde ((x \fiz y^{-1})\fde z^{-1})\right) \right] \fde \left((x\fde y^{-1}) \fiz ((x\fiz y^{-1})\fde z^{-1}) \right).
\end{align*}
We have $(f^{-1}\fde g^{-1}) \fiz ((f^{-1}\fiz g^{-1})\fde h^{-1})
= (f^{-1} \fiz (g^{-1}\fde h^{-1})) \fde (g^{-1}\fiz h^{-1})$
by \eqref{cond-structur3} applied to $\sigma^*$
and suitable $f,g,h \in \Ac$. Filling this identity with $f = x$,
$g^{-1} = x\fde y^{-1}$ and $h^{-1}= (x \fiz y^{-1})\fde z^{-1}$,
we obtain the desired equality.

\bigbreak \emph{Case II.} $(\sigma \times \id) (\id \times
\sigma)(\sigma \times \id) \overset{?}= (\id \times \sigma)(\sigma
\times \id)(\id \times \sigma): \Ac^{\op} \ftimes \Ac \ftimes
\Ac^{\op} \to \Ac^{\op} \ftimes \Ac \ftimes  \Ac^{\op}.$

\bigbreak We fix $(x^{-1},y,z^{-1}) \in \Ac^{\op}\ftimes \Ac
\ftimes \Ac^{\op} $ and compute both sides of the desired
equality:

\begin{align*}
LHS &=
((x^{-1} \fde y) \fde \left( (x^{-1} \fiz y) \fde z^{-1}\right),
(x^{-1}\fde y) \fiz ((x^{-1}\fiz y)\fde z^{-1}),
 (x^{-1}\fiz y) \fiz z^{-1} ),\\
RHS &= (x^{-1} \fde(y \fde z^{-1}),
(x^{-1} \fiz (y\fde z^{-1})) \fde (y\fiz z^{-1}),
(x^{-1} \fiz (y\fde z^{-1})) \fiz (y\fiz z^{-1})).
\end{align*}

The first components of the LHS and the RHS are equal by
\eqref{gstr-aux2}, while the last are equal by \eqref{gstr-aux1}.
It remains to show the equality of the middle components:

\begin{equation}\label{rigid-double-aux2}
(x^{-1}\fde y) \fiz ((x^{-1}\fiz y)\fde z^{-1})
\overset{?}= (x^{-1} \fiz (y\fde z^{-1})) \fde (y\fiz z^{-1}).
\end{equation}
Now
\begin{align*}  LHS \text{ of }\eqref{rigid-double-aux2} &=
\left((x^{-1}\fiz y) \fde y^{-1} \right)^{-1}\fiz ((x^{-1}\fiz
y)\fde z^{-1})
\\ &=  \left\{\left((x^{-1}\fiz y) \fde y^{-1} \right)
\fiz  \left((x^{-1}\fde y) \fde ((x^{-1}\fiz y)\fde
z^{-1})\right)\right\}^{-1}
\\ &=  \left\{\left((x^{-1}\fiz y) \fde y^{-1} \right)
\fiz  \left(x^{-1} \fde (y\fde z^{-1}) \right)\right\}^{-1}
\\ &=  \left\{ \left((x^{-1} \fiz y) \fiz z^{-1}) \right)
\fde  \left(y^{-1} \fiz (y\fde z^{-1}) \right)\right\}^{-1}
\\ &=  \left\{ \left((x^{-1} \fiz (y \fde z^{-1})) \fiz (y \fiz z^{-1})) \right)
\fde  \left(y^{-1} \fiz (y\fde z^{-1}) \right)\right\}^{-1}
\\ &=  \left\{ \left((x^{-1} \fiz (y \fde z^{-1})) \fiz (y \fiz z^{-1})) \right)
\fde  \left(y \fiz z^{-1} \right)^{-1}\right\}^{-1}
\\ &= (x^{-1} \fde (y \fde z^{-1})) \fiz (y \fiz z^{-1}) =
RHS \text{ of }\eqref{rigid-double-aux2}
\end{align*}
Here the first and the second equalities follow from
\eqref{trenzas-notac5}; the third from \eqref{gstr-aux2}; the
fourth from \eqref{cond-structur3} applied to $\sigma^*$ filled
with $f^{-1} = x^{-1} \fiz y$, $g = y$ and $h^{-1}=  y\fde z^{-1}$
as in the first step; the fifth from \eqref{gstr-aux1}; and the
last two from \eqref{trenzas-notac5} again.

\bigbreak \emph{Case III.} $(\sigma \times \id) (\id \times
\sigma) (\sigma\times \id) \overset{?}= (\id \times \sigma)(\sigma
\times \id)(\id \times \sigma): \Ac^{\op} {\,}_\tgt
\hspace{-0.1cm}\times_\sou \Ac^{\op} {\,}_\tgt
\hspace{-0.1cm}\times_\sou  \Ac \to \Ac {\,}_\tgt
\hspace{-0.1cm}\times_\sou \Ac^{\op} {\,}_\tgt
\hspace{-0.1cm}\times_\sou  \Ac^{\op}$.

\bigbreak By Lemma \ref{rigid-dual} (b), $\sigma^{-1}$ is rigid
and hence the equality in case I holds for it. Then case III
follows inverting this equality; indeed
$\left(\sigma^{-1}\right)^* = \left(\sigma^*\right)^{-1}$. \epf

\bigbreak Let now $\Ac$ be a solution. We define  maps
$\sigma^{m,n}: \path_m \Ac \ftimes \path_n \Ac \to \path_n \Ac
\ftimes \path \Ac_m$, $m, n\ge 0$, by
\begin{equation*}
 \sigma^{m,n} =
\left(\sigma_{n, n+1} \dots \sigma_{2,3} \sigma_{1,2} \right)
\left(\sigma_{n+1, n+2} \dots \sigma_{3,4}\sigma_{2,3}  \right)
\dots
\left(\sigma_{n+m-1, n+m} \dots \sigma_{m+1,m+2}\sigma_{m,m+1}  \right),
\end{equation*}
if $m,n>0$; and by
\begin{equation*}
\sigma^{0,n} (\id P, x) = (x, \id Q), \quad
\sigma^{n,0} (x, \id Q) = (\id P,x)
\end{equation*}
if $x\in \path_n \Ac$, $n\ge 0$, $\sou(x) = P$, $\tgt(x) = Q$.
Thus we have an isomorphism of quivers $\psigma: \path \Ac \to
\path \Ac$ by collecting together the maps $\sigma^{m,n}$; and, as
usual, maps $\fde, \fiz, \fd, \fz: \path\Ac {\,}_\tgt
\hspace{-0.1cm}\times_\sou \path\Ac \to\path\Ac$ given by
\eqref{trenzas-notac}.

\begin{lema}\label{rigid-paths}  If $(\Ac, \sigma)$ is a solution
then $(\path \Ac, \psigma)$ is a solution.

If $m, n, p\ge0$ and $(u,v,w) \in \path_m \Ac {\,}_\tgt \hspace{-0.1cm}\times_\sou \path_n \Ac {\,}_\tgt \hspace{-0.1cm}\times_\sou \path_p \Ac$ then
\begin{align}\label{trenza-path1}
u\fde (v\fde w) &= uv\fde w,
\\\label{trenza-path2} u\fiz (v\fiz w) &= u \fiz vw,
\\\label{trenza-path3} u\fde vw &= (u\fde v) \left( (u\fiz v) \fde w \right),
\\\label{trenza-path4} uv\fiz w &= \left(u \fiz (v \fde w)\right)
(v\fiz w).
\end{align}
If $\sigma$ is non-degenerate, so is $\psigma$. \end{lema}

\pf The first claim follows from a well-known equality in the
braid group. Also, $ \sigma^{m+n,p} = (\id\times \sigma^{n,p})
(\sigma^{m,n}\times \id)$, which implies \eqref{trenza-path1} and
\eqref{trenza-path4}, and $ \sigma^{m,n+p} = (\sigma^{m,n}\times
\id) (\id\times \sigma^{n,p})$, which implies \eqref{trenza-path2}
and \eqref{trenza-path3}. Finally, if $u = x_1\dots x_n$ is a path
of length $n$ then the inverse of $u\fde \underline{\;\;}$ is
given by $x_n^{-1} \fde \left(\dots (x_1^{-1}\fde
\underline{\;\;}\right)$ by \eqref{trenza-path1}. Similarly
$\underline{\;\;} \fiz u$ is invertible by \eqref{trenza-path2}.
\epf

The next natural step is to show that $\gstr$ has also the
structure of braided quiver. To state this appropriately we study
in the next subsection the notion of braided groupoid. We come
back to the structure groupoid in Theorem \ref{accion2}.

\bigbreak
\subsection{Braided groupoids}

\

\bigbreak The notion of braided groupoid generalizes the notion of
``braided group" introduced by Takeuchi \cite{tak} to reformulate
results of Lu, Yan and Zhu \cite{lyz1}.

\begin{definition}\label{braidedgpds}
A \emph{braided groupoid} is a collection $(\G, \fde, \fiz)$,
where $\G$ is groupoid and
$$\begin{CD} \G @<\fiz<< \G \ftimes \G @>\fde>>  \G \end{CD}$$
are a left and a right actions,
such that $(\G, \G)$, with $\fde$, $\fiz$,
form a matched pair of groupoids, see
Definition \ref{matchpairgpds}, and
\begin{equation}
\label{braidedgpds-1} fg = (f\fde g)(f \fiz g), \qquad (f,g) \in
\G \ftimes \G.
\end{equation}

A morphism of braided groupoids is a morphism of groupoids that
preserves the actions $\fde$, $\fiz$.
\end{definition}

\bigbreak
Let $\G$ be a braided groupoid. Then the maps $\iota_1, \iota_2:
\G \to \G\bowtie \G$ given by
\begin{equation}\label{iota}
\iota_1(g) = (g, \id \tgt (g)), \qquad \iota_2(g) = (\id \sou (g),
g)
\end{equation}
are morphisms of groupoids.

\bigbreak There is some redundancy in the definition of braided
groupoid, that we study in the next Lemma.

\begin{lema}\label{brgpd-tecnico}
 (a). Let $\G$ be a groupoid  endowed with a left and a right actions
$\fde, \fiz$ such that \eqref{braidedgpds-1} holds. Then
$\G$, $\G$, $\fde$, $\fiz$ form a matched pair of groupoids
(and $\G$ is a braided groupoid).

\bigbreak
(b). Let $\G$, $\G$, $\fde$, $\fiz$ form a matched pair of
groupoids. Then \eqref{braidedgpds-1} holds iff the multiplication
is a morphism of groupoids $\G \bowtie \G \to \G$.

\bigbreak (c). Let $\G$ be a groupoid  endowed with a left action
$\fde: \G \ftimes \G \to \G$. Let $\fiz: \G \ftimes \G \to \G$ be
given by \eqref{braidedgpds-1}, \emph{i. e.} $x\fiz y = (x\fde
y)^{-1}xy$, $(x,y)\in \G\ftimes\G$. Then $\fiz$ is a right action
if and only if \eqref{mp-3} holds. If this is the case, then $\G$,
$\G$, $\fde$, $\fiz$ form a matched pair of groupoids, and $\G$ is
a braided groupoid.

\bigbreak (d). Let $\G$ be a groupoid  endowed with a right action
$\fiz: \G \ftimes \G \to \G$. Let $\fde: \G \ftimes \G \to \G$ be
given by \eqref{braidedgpds-1}, \emph{i. e.} $x\fde y = xy(x\fiz
y)^{-1}$, $(x,y)\in \G\ftimes\G$. Then $\fde$ is a left action if
and only if \eqref{mp-4} holds. If this is the case, then $\G$,
$\G$, $\fde$, $\fiz$ form a matched pair of groupoids, and $\G$ is
a braided groupoid.
\end{lema}

\pf  (a).  It is clear that \eqref{mp-0.7} holds. We check \eqref{mp-3}:
\begin{align*}
(f\fde (gh))(f\fiz (gh)) &= fgh = (f\fde g)(f \fiz g)h
\\ &= (f\fde g)[(f \fiz g) \fde h][(f \fiz g) \fiz h]
= (f\fde g)[(f \fiz g) \fde h](f\fiz (gh)).
\end{align*}
The proof of \eqref{mp-4} is similar; thus (a) is valid.
The proofs of (b), (c) and (d) are straightforward. \epf

The antipode $x \mapsto x^{-1}$ induces an isomorphism of quivers
$\G \to \G^{\op}$. We check that the identities discussed in Lemma
\ref{rigid-dual} are valid with respect to the antipode.

\begin{lema}\label{antipode2}
Let $(\G, \fde, \fiz)$ be a braided groupoid.

(a). The identities \eqref{trenzas-notac3},
\eqref{trenzas-notac4}, \eqref{trenzas-notac5} hold in $\G$ with
respect to the antipode.

(b). $(\G, \fd, \fz)$ is also a braided groupoid.
\end{lema}

\pf (a). The validity of \eqref{trenzas-notac3} is clear since
$\fde$, $\fiz$ are actions. If $(x,g^{-1}) \in \G \ftimes \G$ then
$$
\id \sou(x) = x\fde \id \tgt(x) = x\fde (g^{-1}g) = (x\fde g^{-1})((x\fiz g^{-1}) \fde g),
$$
similarly for the other identity and \eqref{trenzas-notac5}
follows. Let now $\widetilde\sigma:  \Ac \ftimes \Ac \to  \Ac
\ftimes \Ac$ be given by $\widetilde\sigma (x,y) =
\left((y^{-1}\fiz x^{-1})^{-1}, (y^{-1}\fde x^{-1})^{-1}\right)$,
$(x,y) \in \G \ftimes \G$. We compute:
\begin{align*}
\sigma \widetilde\sigma (x,y) &= \sigma \left((y^{-1}\fiz
x^{-1})^{-1}, (y^{-1}\fde x^{-1})^{-1}\right)
\\ &= \left((y^{-1}\fiz x^{-1})^{-1}
\fde(y^{-1}\fde x^{-1})^{-1}, (y^{-1}\fiz x^{-1})^{-1} \fiz
(y^{-1}\fde x^{-1})^{-1}\right).
\end{align*}
Now
\begin{align*}
(y^{-1}\fiz x^{-1})^{-1} \fde(y^{-1}\fde x^{-1})^{-1} &=
(y\fiz (y^{-1} \fde x)) \fde ((y^{-1}\fiz x^{-1}) \fde x)
\\ &= \left((y\fiz (y^{-1} \fde x))(y^{-1}\fiz x^{-1}) \right) \fde x
\\ &= (yy^{-1}\fiz x^{-1})  \fde x = x;
\\ (y^{-1}\fiz x^{-1})^{-1} \fiz (y^{-1}\fde x^{-1})^{-1}
&= (y\fiz (y^{-1} \fde x)) \fiz ((y^{-1}\fiz x^{-1}) \fde x)
\\ &= y\fiz \left( (y^{-1} \fde x) \left((y^{-1}\fiz x^{-1}) \fde x \right) \right)
\\ &= y\fiz \left( y^{-1} \fde x x^{-1}  \right) = y.
\end{align*}
Similarly, $\widetilde\sigma \sigma = \id$. Thus,
\eqref{trenzas-notac4} holds.

\bigbreak (b). We have $(x\fd y)(x\fz y) = \left((x\fd y)\fde
(x\fz y)\right) \left((x\fd y)\fiz (x\fz y)\right) = xy$. We claim
that $\fd$ is a left action. Let $(x,y,z) \in \G \ftimes \G\ftimes
\G$. Then $$x \fd (y\fd z) = x \fd (y^{-1}\fiz z^{-1})^{-1} =
\left((y^{-1}\fiz z^{-1})\fiz x^{-1}\right)^{-1} =
\left(y^{-1}\fiz z^{-1}x^{-1}\right)^{-1} = xy \fd z.$$ Similarly,
$\fz$ is a right action. By Lemma \ref{brgpd-tecnico} (a), $(\G,
\fd, \fz)$ is a braided groupoid. \epf

\bigbreak
\subsection{Braided groupoids are braided quivers}

\

\bigbreak We next justify the name of ``braided groupoids".

\begin{lema} Let $\G$ be a braided groupoid.
Let $\sigma:  \G \ftimes \G \to  \G
\ftimes \G$ be the map
\begin{equation}
\label{braidedgpds-2} \sigma (f,g) = (f\fde g, f \fiz g), \qquad
(f,g) \in \G {\,}_\tgt \hspace{-0.1cm}\times_\sou \G.
\end{equation}
Then $\sigma$ is a non-degenerate solution of the
quiver-theoretical braid equation \eqref{qtBE}.
\end{lema}

\pf We fix $(f,g,h) \in  \G \ftimes
\G {\,}_\tgt \hspace{-0.1cm}\times_\sou  \G$ and compute:
\begin{align*}
(\sigma \times \id) (\id \times \sigma)(\sigma \times \id)  (f,g,h) &=
((fg) \fde h, (f\fde g) \fiz ((f\fiz g)\fde h), f\fiz (gh)),\\
(\id \times \sigma)(\sigma \times \id)(\id \times \sigma) (f,g,h) &=
((fg) \fde h,  (f \fiz (g\fde h)) \fde (g\fiz h), f\fiz (gh)).
\end{align*}
Now we apply several times \eqref{braidedgpds-1} and compute:
\begin{align*}
fgh &= (f\fde g)(f \fiz g)h
\\ &= (f\fde g)[(f \fiz g) \fde h][(f \fiz g) \fiz h]
\\ &= [(f\fde g) \fde ((f\fiz g)\fde h)] [(f\fde g) \fiz ((f\fiz g)\fde h)] [f\fiz (gh)]
\\ &= [(fg) \fde h] [(f\fde g) \fiz ((f\fiz g)\fde h)] [f\fiz (gh)]
\end{align*}
and also
\begin{align*}
fgh &= f(g\fde h)(g\fiz h)
\\ &= [f\fde(g\fde h)][f\fiz (g\fde h)] (g\fiz h)
\\ &= [(fg) \fde h]  [(f\fiz (g\fde h)) \fde (g\fiz h)] [(f\fiz (g\fde h)) \fiz (g\fiz h)]
\\ &=  [(fg) \fde h] [(f\fiz (g\fde h)) \fde (g\fiz h)] [f\fiz (gh)].
\end{align*}
Hence $(f\fde g) \fiz ((f\fiz g)\fde h) = (f\fiz (g\fde h)) \fde
(g\fiz h)$, and $\sigma$ is a solution of \eqref{qtBE}. Since
$\fde$, $\fiz$ are actions, $\sigma$ is non-degenerate. \epf

\bigbreak
Let us say that a sub-quiver $\Ac$ of a braided groupoid $\G$ is
\emph{invariant} if $\Ac \fde \Ac \subset \Ac$, $\Ac \fiz \Ac
\subset \Ac$, $\Ac^{-1} \fde \Ac \subset
\Ac$ and $\Ac \fiz \Ac^{-1} \subset \Ac$.

\begin{cor} Let $\Ac$ be an invariant sub-quiver of a braided groupoid $\G$.
Then $(\Ac, \sigma_{\vert \Ac {\,}_\tgt \hspace{-0.1cm}
\times_\sou \Ac})$ is a non-degenerate braided quiver.
\qed\end{cor}

Thus, braided groupoids and their invariant sub-quivers are
naturally braided quivers. See Remark \ref{faithful} below.

\bigbreak In the papers \cite{ess, lyz1, s}, braided structures on
groups were described through suitable 1-cocycles. We show now
that this description goes over also to groupoids but with group
bundles as recipients of the 1-cocycles.

\begin{definition} A \emph{1-cocycle groupoid datum} is
a triple $(\G, \N, \pi)$ where $\G$ is a groupoid over $\Pc$, $\N$
is a group bundle over $\Pc$, provided with a right action $\fiz:
\N {\,}_{p} \hspace{-0.1cm}\times_\sou \G \to \N$ by group bundle
automorphisms; and $\pi: \G \to \N$ is a bijective 1-cocycle,
\emph{i.e.} it is a bijection with $p\pi = {\sou}$ and
\begin{equation}\label{cociclo}
\pi(fg) = (\pi(f) \fiz g) \pi(g), \qquad (f,g) \in \G
{\,}_{\tgt} \hspace{-0.1cm} \times_\sou \G.
\end{equation}
\end{definition}

\begin{theorem}\label{coc-gpd-datum} Let $\G$ be a groupoid.
There is a bijective correspondence between
\begin{enumerate}
\item[(a)] Structures of braided groupoid $(\G, \sigma)$.

\bigbreak
\item[(b)] 1-cocycle groupoid data $(\G, \N, \pi)$.
\end{enumerate}

\bigbreak
In this correspondence, $\G \bowtie \G \simeq \G \ltimes \N$, and
\begin{equation}\label{trenza}
\sigma(f, g) = \left(fg\left(\pi^{-1}\left(\pi(f) \fiz
g \right)\right)^{-1}, \pi^{-1}\left(\pi(f) \fiz g \right)\right), \qquad (f,g) \in \G
{\,}_{\tgt} \hspace{-0.1cm} \times_\sou \G.
\end{equation}

\end{theorem}

\pf (a) $\implies$ (b). Let $\N$ be the kernel of the
multiplication map $\G \bowtie \G \to \G$; identify $\G$ with
$\G_1 :=$ the  image of $\iota_1$, \emph{cf.} \eqref{iota}; let
$\fiz$ be the restriction of the adjoint action. Finally, let
$\pi: \G \to \N$ be defined by $\pi(g) = (g^{-1}, g)$, $g\in \G$.
Then $(\G, \N, \pi)$ is a 1-cocycle groupoid datum. Indeed,
\begin{align*}
 (\pi(f) \fiz g) \pi(g) &= (g^{-1}, \id {\sou} (g))
 (f^{-1}, f)(g, \id \tgt (g))(g^{-1}, g)
\\ &= (g^{-1} (\id {\sou} (g)\fde f^{-1}), (\id {\sou} (g)\fiz f^{-1}) f)
(g (\id \tgt (g)\fde g^{-1}), (\id \tgt (g)\fiz g^{-1}) g^{-1})
\\ &= (g^{-1}f^{-1},  f) (\id {\sou} (g), g) =
(g^{-1}f^{-1}(f \fde \id {\sou} (g)), (f \fiz \id {\sou} (g))g)
\\ &= (g^{-1}f^{-1},  fg) = \pi(f,g).
\end{align*}

(b) $\implies$ (a). Consider the map $\psi: \G {\,}_\tgt
\hspace{-0.1cm}\times_\sou \G \to \G \ltimes \N$, $\psi(f,g) =
(fg, \pi(g))$, $(f,g) \in \G {\,}_\tgt \hspace{-0.1cm}\times_\sou
\G $. Let $\G_1$, resp. $\G_2$, be the image of $\iota_1$, resp.
$\iota_2$, as in \eqref{iota}. Then
$$
\psi(\G_1) = \{(g, \id \tgt (g)): g\in \G\}, \qquad
\psi(\G_2) = \{(g, \pi(g)): g\in \G\}
$$
are subgroupoids of $\G\ltimes \N$ isomorphic to $\G$, and they
form an exact factorization of $\G\ltimes \N$. Transporting the
structure back to $\G {\,}_\tgt \hspace{-0.1cm}\times_\sou \G$ via
$\psi$, we have a groupoid structure on this, which is isomorphic
to $\G_1\bowtie \G_2 \simeq \G\bowtie \G$. A straightforward computation
shows that the induced actions
$\begin{CD} \G @<\fiz<< \G \ftimes \G @>\fde>>  \G \end{CD}$
are explicitly given by
\begin{equation}\label{acciones}
f \fde g = fg\left(\pi^{-1}\left(\pi(f) \fiz g \right)\right)^{-1},
\qquad f \fiz g = \pi^{-1}\left(\pi(f) \fiz g \right), \qquad (f,g) \in \G
{\,}_{\tgt} \hspace{-0.1cm} \times_\sou \G.
\end{equation}
Hence \eqref{braidedgpds-1} holds, and the corresponding solution is
given by \eqref{trenza}.

\bigbreak Finally, it is not difficult to see that these
constructions  are inverse of each other. \epf

It follows from \eqref{acciones} that
\begin{equation}\label{dentro}
\pi(f\fiz g) = \pi(f) \fiz g, \qquad (f,g) \in \G
{\,}_{\tgt} \hspace{-0.1cm} \times_\sou \G.
\end{equation}

\bigbreak

A \emph{symmetric groupoid} is a braided groupoid such that the
corresponding $\sigma$ is a symmetry.
The next characterization of symmetric groupoids generalizes results from \cite{ess, lyz1}.

\begin{proposition} A braided groupoid $\G$ is symmetric if and only if
the corresponding $\n$ is abelian.
\end{proposition}

\pf Let $(f,g) \in \G {\,}_\tgt \hspace{-0.1cm}\times_\sou \G$. Then
\begin{align*}
\pi(f\fiz g )  \pi(g) &= (\pi(f) \fiz g) \pi(g) = \pi(fg)
= \pi \left((f\fde g)(f\fiz g)\right)
\\  & =  \left(\pi(f\fde g) \fiz(f\fiz g)\right) \pi(f\fiz g)
= \pi\left((f\fde g) \fiz(f\fiz g)\right) \pi(f\fiz g),
\end{align*}
by \eqref{dentro}, \eqref{cociclo} and \eqref{braidedgpds-1}. This
says that $\n$ is abelian if and only if $g = (f\fde g) \fiz(f\fiz
g)$ for any $(f,g) \in \G {\,}_\tgt \hspace{-0.1cm}\times_\sou
\G$. But $\sigma^2(f,g) = \sigma(f \fde g, f\fiz g) = \left((f
\fde g)\fde (f\fiz g), (f \fde g)\fiz (f\fiz g)\right)$. Hence, if
$\G$ is symmetric then $\n$ is abelian. If $\n$ is abelian, then
$fg = \left((f \fde g)\fde (f\fiz g)\right) \left((f \fde g)\fiz
(f\fiz g)\right)= \left((f \fde g)\fde (f\fiz g)\right)g$, thus
$\G$ is symmetric. \epf

\bigbreak We define next the subgroup bundles $\Gamma_r :=$ kernel
of $\fiz$, $\Gamma_l :=$  kernel of $\fde$ and $\Gamma :=\Gamma_r
\cap \Gamma_l$ of the braided groupoid $\G$. Hence

\begin{equation}\label{gamadef}
\Gamma = \{v \in \G: \, v\fde w = w, \quad z\fiz v = z, \quad
\forall w,z \in \G, \quad \sou(w) = \tgt(v) = \sou(v)  = \tgt(z) \}.
\end{equation}

By \eqref{braidedgpds-1}, we also have
\begin{align*}
\Gamma_l &= \{m \in \G: m\fiz y = y^{-1}my, \quad \forall y \in \G,
\quad \sou(m) = \tgt(m) = \sou(y) \},
\\ \Gamma_r &= \{n \in \G: \,\,\, x\fde n = xnx^{-1}, \,\quad
\forall x \in \G, \;\quad \tgt(x) = \sou(n) \;=\tgt(n) \}.
\end{align*}

\begin{lema}\label{gama} Let $\G$ be a braided groupoid.

\begin{enumerate}
\item[(a)] The wide subgroup bundle $\Gamma$ defined above is  abelian and normal.

\bigbreak
\item[(b)] If $\Lambda \subset \Gamma$ is a normal subgroup bundle of $\G$
then $\G / \Lambda$ is a braided groupoid, with braiding inherited
from $\G$ via the canonical projection.
\end{enumerate}
\end{lema}

\pf (a). It follows at once from \eqref{braidedgpds-1} that $\Gamma$ is
abelian, and it is clearly normal.

(b). One checks that $\Gamma_l {\,}_\tgt
\hspace{-0.1cm}\times_\sou \Pc$ and $\Pc {\,}_\tgt
\hspace{-0.1cm}\times_\sou \Gamma_r$ are normal subgroup bundles
of $\G \bowtie \G$; hence $\Gamma {\,}_\tgt
\hspace{-0.1cm}\times_\sou \Gamma$, and \emph{a fortiori} $\Lambda
{\,}_\tgt \hspace{-0.1cm}\times_\sou \Lambda$, are normal subgroup
bundles of $\G \bowtie \G$. Then the quotient $\G \bowtie \G /
\Lambda {\,}_\tgt \hspace{-0.1cm}\times_\sou \Lambda$ factors as a
product of the wide subgroupoids $\G / \Lambda {\,}_\tgt
\hspace{-0.1cm}\times_\sou \Pc$ and $\Pc \ftimes \G / \Lambda$;
this factorization induces a matched pair structure, and hence a
structure of braided groupoid, on $\G / \Lambda$. \epf

We close this subsection with an application of Theorem
\ref{coc-gpd-datum}; this is a generalization of the examples in
\cite[Section 3]{lyz1}, \cite{wx}.

\bigbreak Let $(\Vc, \Hc)$ be a matched pair of groupoids. Recall
that the restricted product of $\Vc$ and $\Hc$ is the groupoid
$\Vc \prode \Hc := \{(g, x) \in  \Vc \times \Hc: \sou(g) =
\sou(x), \quad \tgt(g) = \tgt(x)\}$, with component-wise product
\cite{AA}.

\begin{proposition} Let $\G := \Vc \prode \Hc$,
$\N := (\Vc \bowtie \Hc)^{\bundle}$; let $\pi: \G \to \N$ and
$\fiz: \N {\,}_{p} \hspace{-0.1cm}\times_\sou \G \to \N$  be given
by
\begin{align} \label{wx1}
\pi(g,x) &= g^{-1}x,
\\ \label{wx2}
d \fiz (g,x) &= g^{-1}dg,
\end{align}
$(g,x)\in \Vc \prode \Hc$, $d\in (\Vc \bowtie \Hc)^{\bundle}$.
Then $(\G, \N, \pi)$ is a 1-cocycle groupoid datum. Thus, $(\Vc
\prode \Hc, \sigma)$, where
\begin{equation}\label{trenza-mp}
\sigma\left((g,x), (h,y)\right)
= \left(\left(x\fde h, xy(x\fiz h)^{-1}\right),\left((x\fde h)^{-1}gh,
x\fiz h \right)\right),
\end{equation}
$(g,x), (h,y) \in \Vc \prode \Hc$, $\tgt(x) = \sou (h)$,  is a
braided groupoid.
\end{proposition}

\pf A straightforward verification shows that $(\G, \N, \pi)$ is a
1-cocycle groupoid datum. The explicit formula \eqref{trenza-mp}
follows from \eqref{trenza}, once we show that the actions $\fde$,
$\fiz$ of $\Vc \prode \Hc$ are given by
\begin{equation}\label{acciones-mp}
(g,x) \fde (h,y) = \left(x\fde h, xy(x\fiz h)^{-1}\right),
\qquad (g,x) \fiz (h,y) =\left((x\fde h)^{-1}gh, x\fiz h \right),
\end{equation}
$(g,x), (h,y) \in \Vc \prode \Hc$, $\tgt(x) = \sou (h)$.
By \eqref{acciones}, we have
\begin{align*} (g,x) \fiz (h,y) &= \pi^{-1}\left(\pi(g,x) \fiz (h,y) \right) =
\pi^{-1}\left(g^{-1}x \fiz (h,y) \right) =
\pi^{-1}\left(h^{-1}g^{-1}x h \right)
\\ & = \pi^{-1}\left(h^{-1}g^{-1}(x\fde h) (x\fiz h) \right)
= \left( (x\fde h)^{-1} gh, x\fiz h\right).
\end{align*}
This shows the second equality in \eqref{acciones-mp}, and implies
the first by \eqref{acciones}. \epf

\section{Characterizations of braided quivers}

\bigbreak
\subsection{LYZ-pairs}

\

\bigbreak We begin by a categorical way of constructing braided
quivers. Let $(\Vc, \Hc)$ be a matched pair of groupoids and
recall the definition of representation of $(\Vc, \Hc)$ in
Subsection \ref{mpg}.

\begin{definition}\label{lyzpairs} (\cite{AA}, inspired in \cite{lyz3, tak}).
Let $\kappa: \Vc \to \Hc$ be a morphism of groupoids. We shall say
that $\kappa$ is a \emph{rotation} if
\begin{equation}\label{lyz0}
y\kappa(g) = \kappa(y\fde g)(y\fiz g) \qquad \text{ for all } g\in
\Vc, \, y \in \Hc, \, b(g) = l(y).
\end{equation}
Let $\lambda: \Vc \bowtie \Hc \to \Hc$, given by $\lambda(g, x) =
\kappa(g)x$. Then $\kappa$ is a rotation if and only if $\lambda$
is a morphism of groupoids, see \cite{AA}.

\bigbreak A \emph{LYZ-pair}\footnote{LYZ pairs are called ``matched pairs of rotations" in \cite{AA}.} is a pair $(\lyzu,
\lyzd)$ of rotations $\Vc \to \Hc$ such that
\begin{align}
\label{lyz1} \eta(g) \fde f = gf \bigl(\xi(f)^{-1} \fde
g^{-1}\bigr)\,.
\end{align}
for every $f$ and $g$ in $\Vc$ with $b(g)=t(f)$.
\end{definition}

\begin{theorem}\label{clasifbrstr} \cite{AA}. Structures of
braided category on $\Rep(\Vc, \Hc)$ are parameterized by LYZ-pairs.
If $\Ac, \Bc$ are representations of $(\Vc,
\Hc)$ and $(\lyzu, \lyzd)$ is a LYZ-pair, then
the induced braiding $\sigma_{\Ac, \Bc}: \Ac
{\,}_\tgt\hspace{-0.1cm} \times_\sou \Bc \to \Bc
{\,}_\tgt\hspace{-0.1cm} \times_\sou \Ac$ is given by
\begin{equation}\label{trenza23}
\sigma_{\Ac, \Bc} (a, b)  = \Bigl(\lyzd(\vert a\vert) \fde b,
\left( \lyzu(\vert b\vert)^{-1} \fiz \vert a \vert^{-1}\right)
\fde a \Bigr), \quad (a, b)\in \Ac {\,}_\tgt\hspace{-0.1cm}
\times_\sou \Bc.
\end{equation}
\qed
\end{theorem}

It can be easily shown that
\begin{equation}\label{trenza24}
\sigma^{-1}_{\Bc, \Ac} (a, b)  = \Bigl(\lyzu(\vert a\vert) \fde b,
\left( \lyzd(\vert b\vert)^{-1} \fiz \vert a \vert^{-1}\right)
\fde a \Bigr), \quad (a, b)\in \Ac {\,}_\tgt\hspace{-0.1cm}
\times_\sou \Bc.
\end{equation}

\begin{cor}\label{exabrq} Let $(\Vc, \Hc)$ be a matched pair of
groupoids and let $(\lyzu, \lyzd)$ be a LYZ-pair.
If $\Ac$ is a representation of $(\Vc, \Hc)$ then $(\Ac,
\sigma_{\Ac, \Ac})$ is a braided quiver.
\end{cor}
\pf
The statement follows from a well-known
general result in braided categories.
\epf

As we shall see in the next Subsection, any non-degenerate braided
quiver arises in this way. To this end we shall need the following
result that generalizes \cite[Prop. 5.1]{tak}.

\begin{theorem}\label{lematak2}
Let $\G$ be a braided groupoid and let $\G\bowtie \G$ be the
corresponding diagonal groupoid with respect to $\fde$, $\fiz$.

Let $\fdd: \G \bowtie \G \ftimes \G  \to \G$, $\fzz: \G \bowtie \G
\ftimes \G  \to \G \bowtie \G$,  and $\lyzuu, \lyzdd: \G \to \G
\bowtie \G$ be given by
\begin{align}\label{accionbt1}
(g,h) \fdd f &= g\fd(h \fde f),
\\ (g,h) \fzz f &= (g \fz(h \fde f), h\fiz f),
\\ \lyzuu(f) &= (f, \id \tgt(f)),
\\ \lyzdd(f) &= (\id \sou(f), f),
\end{align}
$(g,h)\in \G \bowtie \G$,  $f \in\G$, $\tgt(h) = \sou(f)$. Then
$(\G, \G \bowtie \G)$, with $\fdd$, $\fzz$, is a matched pair of
groupoids, and $(\lyzuu, \lyzdd)$ is a LYZ-pair
for it.
\end{theorem}

\pf  We first check that $\fdd$ is a left action. If $(g,h, f, k,
\ell)\in \G \ftimes \G\ftimes \G\ftimes \G\ftimes \G$ then
\begin{align*}
(g,h)\fdd \left( (f,k) \fdd \ell \right) &= (g,h)\fdd \left(
f\fd(k \fde \ell) \right) = g\fd \left[h \fde \left( f\fd(k \fde
\ell) \right)\right];
\\  (g,h) (f,k) \fdd \ell  &= \left(g(h\fde f), (h\fiz f)k \right) \fdd \ell
= \left[g(h\fde f)\right] \fd \left[\left((h\fiz f)k \right)\fde \ell \right]
\\  &= g\fd \left\{(h\fde f) \fd \left[(h\fiz f)\fde (k\fde \ell) \right]\right\}
\end{align*}
Thus $\fdd$ is a left action by \eqref{tecnica1}.
We next check that $\fzz$ is a right action. If $g, h f, k$ are as above, then
\begin{align*}
\left((g,h) \fzz f\right) \fzz k &= (g \fz(h \fde f), h\fiz f)\fzz k
\\ &= \left( \left((g \fz(h \fde f)\right) \fz \left((h\fiz f) \fde k\right), (h\fiz f) \fiz k\right)
\\  &= \left(g \fz \left((h \fde f)\left((h\fiz f) \fde k\right)\right), h\fiz fk\right)
\\  &=
\left(g \fz \left(h \fde f k\right), h\fiz fk\right)
= (g,h) \fzz fk.
\end{align*}

We next check the compatibility conditions, \eqref{mp-0.7} being
clear. We verify \eqref{mp-3}; let $g,h,f,k$ as before. Then
\begin{align*}
\left((g,h) \fdd f \right) \left[\left((g,h) \fzz f \right) \fdd k\right]
&= \left[g\fd(h \fde f)\right]\left[(g \fz(h \fde f), h\fiz f)\fdd k\right]
\\  &= \left[g\fd(h \fde f)\right] \left\{\left[(g \fz(h \fde f)\right] \fd \left[(h\fiz f)\fde k \right]\right\}
\\  &= g \fd \left[(h \fde f)((h \fiz f) \fde k)\right]
\\  &= g \fd (h \fde fk)  = (g,h) \fdd fk.
\end{align*}
If $g,h,f,k,\ell$ are as above, then the left-hand side of \eqref{mp-4} is
\begin{align*}
(g,h) (f, k) \fzz \ell &= \left(g(h\fde f), (h\fiz f) k \right)
\fzz \ell
\\  & = \left(\left(g(h\fde f)\right) \fz \left[\left((h\fiz f) k\right) \fde \ell \right],
\left((h\fiz f) k\right) \fiz \ell\right),
\end{align*}
whose first component is
\begin{align*}
\left(g(h\fde f)\right) \fz \left[\left((h\fiz f) k\right) \fde \ell \right]
&= \left\{g \fz  \left[(h\fde f)\fd \left(\left((h\fiz f) k\right) \fde \ell \right) \right]\right\}
\left\{(h\fde f) \fz \left[\left((h\fiz f) k\right) \fde \ell \right] \right\}
\\  & = \left\{g \fz  \left[(h\fde f)\fd \left[(h\fiz f)\fde (k
\fde \ell) \right] \right]\right\}
\\ & \qquad\qquad\qquad
\left\{(h\fde f) \fz \left[(h\fiz f)\fde (k
\fde \ell) \right] \right\}
\end{align*}

\bigbreak
On the other hand, the first component of the right-hand side of \eqref{mp-4}:
\begin{align*}
\left((g,h) \fzz \left((f, k) \fdd \ell\right)\right)
\left((f, k) \fzz \ell\right) &=
\left(g \fz \left[h \fde \left[f\fd (k \fde \ell)\right]\right],
h \fiz \left[f\fd (k \fde \ell)\right] \right)
(f \fz (k \fde \ell), k\fiz \ell)
\end{align*}
is
\begin{align*}
& \left\{g \fz \left[h \fde \left[f\fd (k \fde \ell)\right]\right]\right \}
\left\{\left[ h \fiz \left[f\fd (k \fde \ell)\right]\right]
\fde \left[f \fz (k \fde \ell)\right]\right\}
\end{align*}
and we have equality with the former because of \eqref{tecnica1} and
\eqref{tecnica2}.
Finally, the second component of the right-hand side of \eqref{mp-4}
is:
\begin{align*}
\left\{\left[ h \fiz \left[f\fd (k \fde \ell)\right]\right]
\fiz \left[f \fz (k \fde \ell)\right]\right\} (k\fiz \ell)
& = \left\{ h \fiz \left[\left[f\fd (k \fde \ell)\right]
 \left[f \fz (k \fde \ell)\right]\right]\right\} (k\fiz \ell)
\\ & = \left\{ h \fiz \left[f (k \fde \ell)\right]\right\} (k\fiz \ell)
\\ & = \left\{ (h \fiz f) \fiz (k \fde \ell)\right\} (k\fiz \ell)
\\ & = \left( (h \fiz f) k \right) \fiz \ell.
\end{align*}
Thus, we have shown the validity of \eqref{mp-4}.

\bigbreak We claim now that the maps $\lyzuu$ and $\lyzdd$ are
rotations. If $g,h,f$ are as above then
\begin{align*}
(g,h) \lyzuu(f) &= (g,h) \left(f, \id \tgt(f)\right)
=\left(g(h\fde f), (h\fiz f)\right);
\\ \lyzuu\left(((g,h)\fdd f\right)\left((g,h)\fzz f\right)
&= \left(g\fd(h \fde f), \id \tgt(g\fd(h \fde f))\right)\left(g \fz(h \fde f), h\fiz f\right)
\\ &= \left(\left(g\fd(h \fde f)\right)\left(g \fz(h \fde f)\right),
h\fiz f\right)=\left(g(h\fde f), (h\fiz f)\right);
\end{align*}
\begin{align*}
\\(g,h) \lyzdd(f) &= (g,h) \left(\id \sou(f), f\right)
=\left(g, hf\right);
\\ \lyzdd\left(((g,h) \right.& \left.
\fdd f\right)\left((g,h)\fzz f\right)
= \left(\id \sou(g\fd(h \fde f)), g\fd(h \fde f)\right)\left(g \fz(h \fde f), h\fiz f\right)
\\ &= \left(\left(g\fd(h \fde f)\right) \fde \left(g \fz(h \fde f)\right), \left[\left(g\fd(h \fde f)\right) \fiz \left(g \fz(h \fde f)\right)\right](h\fiz f)\right)
\\ &= \left(g, (h \fde f)(h\fiz f)\right),
\end{align*}
where in the last equality we have used \eqref{trenzas-notac1.5}.

\bigbreak We next verify the condition \eqref{lyz1}.
\begin{align*}
gf \left(\lyzuu(f)^{-1} \fde g^{-1} \right) &=
gf \left( (f^{-1}, \id \sou(f)) \fdd g^{-1} \right)
= gf \left( f^{-1} \fd g^{-1} \right)
\\ &= (g\fde f)(g\fiz f) \left( f^{-1} \fd g^{-1} \right)
= g \fde f;
\\\lyzdd(g) \fdd f &= (\id \sou(g), g) \fdd f = g \fde f.
\end{align*}
We have proved that $(\lyzuu, \lyzdd)$ is a LYZ-pair
for $(\G, \G \bowtie \G)$. \epf

\begin{obs}\label{lematak3}
Let $\G$ be a braided groupoid. Let $\fddd: \G \bowtie \G \ftimes
\G  \to \G$, $\fzzz: \G \bowtie \G \ftimes \G  \to \G \bowtie \G$,
be given by
\begin{align*}
(g,h) \fddd f &= gh \fde f,
\\ (g,h) \fzzz f &= (g \fiz(h \fde f), h\fiz f),
\end{align*}
$(g,h)\in \G \bowtie \G$,  $f \in\G$, $\tgt(h) = \sou(f)$. Then
$(\G, \G \bowtie \G)$, with $\fddd$, $\fzzz$, is a matched pair of
groupoids. We shall not need this result in the sequel, thus we
leave the proof to the reader.
\end{obs}

\begin{obs}\label{lematak1}
Let $(\Vc, \Hc)$ be a matched pair of groupoids
and let $(\lyzu, \lyzd)$ be a LYZ-pair. Then
there is a structure of braided groupoid on $\Vc$ such that $(\Vc,
\Vc \bowtie \Vc)$, with $(\lyzuu, \lyzdd)$, ``covers" $(\Vc, \Hc)$
with $(\lyzu, \lyzd)$. Compare with \cite[Section 5]{tak}. We
shall not need this result in the sequel, so we do not discuss it
in detail.
\end{obs}

\begin{lema}\label{acciones-compatibles}
Let $\G$ be a braided groupoid and let $\Ac$ be a quiver.

\bigbreak
(i). There is a bijective correspondence between
\begin{enumerate}
\item[(a)] Left actions $\fdd$ of $\G \bowtie \G$ on $\Ac$.

\bigbreak
\item[(b)] Pairs $(\fde, \fd)$ of left actions of $\G$ on $\Ac$
such that
\begin{align}\label{acciont-quiver1}
h \fde \left(f \fd y \right) &= \left(h\fde f\right) \fd \left[(h\fiz f)\fde y \right],
\end{align}
\end{enumerate}

Assume that $\Ac$ is a representation of $(\G, \G \bowtie \G)$.

\bigbreak
(ii). The left action $\fd$ is compatible with the definition \eqref{trenzas-notac1}.

\bigbreak
(iii). The map $\vert \ \vert$ is a morphism of braided quivers.
Therefore it preserves $\fde, \fiz, \fd, \fz$.

\bigbreak
(iv). The left action $\fdd$ of $\G \bowtie \G$ on $\Ac$ induces a
left action of $\G$ on $\Ac \ftimes \Ac$ by
\begin{equation}\label{accion-diagonal}
g \fdd (x,y) := (g\fde x, (g\fiz\vert x\vert) \fd y),
\end{equation}
$g\in \G$, $(x,y) \in \Ac \ftimes \Ac$, where $\fde$, $\fd$ are as
in (b).
\end{lema}

We shall denote by $\nabla: \G \to \aut (\Ac \ftimes \Ac)$ the map
induced by the action \eqref{accion-diagonal}.

\pf (i).
(b) $\implies$ (a). The correspondence is given by
\begin{equation}\label{acciont-quiver0}
(g,h) \fdd x = g\fd(h \fde x),
\end{equation}
$(g,h)\in \G \bowtie \G$,  $x \in\Ac$, $\tgt(h) = \sou(x)$. We
check that $\fdd$ is a left action. If $(g,h), (f, k)\in \gstr
\bowtie \gstr$, $x \in\Ac$, $\tgt(k) = \sou(x)$, then
\begin{align*}
(g,h)\fdd \left( (f,k) \fdd x\right) &= (g,h)\fdd \left(f \fd (k \fde x) \right)
= g\fd \left[h \fde \left(f \fd (k \fde x) \right)\right];
\\  (g,h) (f,k) \fde x  &= \left(g(h\fde f),(h\fiz f)k \right) \fde x
= \left(g(h\fde f)\right) \fd \left[\left((h\fiz f)k  \right)\fde x \right]
\\  &=  g \fd \left\{\left(h\fde f\right) \fd \left[(h\fiz f)\fde (k\fde x) \right] \right\}.
\end{align*}
Letting $y = k\fde x$, we are reduced to \eqref{acciont-quiver1}.
(a) $\implies$ (b) is similar.

\bigbreak
(ii). We have $a\fd b= \lyzuu(a) \fdd b = \vert a \vert \fd b$,
the first equality by \eqref{trenza24}.

\bigbreak
(iii). We have $\sigma(a,b) =  \Bigl(\lyzdd(\vert a\vert) \fdd b,
\left( \lyzuu(\vert b\vert)^{-1} \fzz \vert a \vert^{-1}\right)
\fdd a \Bigr) = \Bigl(\vert a\vert \fde b, \left(\vert b\vert^{-1}\fz \vert a\vert^{-1}\right) \fd a \Bigr)$, hence $\vert \ \vert$ is a morphism of braided quivers.

\bigbreak
(iv). Straightforward using \eqref{compcond}. \epf

\bigbreak
\subsection{Groupoid-theoretical characterization of braided quivers}

\

\bigbreak We first state a characterization of the structure
groupoid by a universal property, a generalization of \cite[Th.
9]{lyz1}. Let $(\Ac, \sigma)$ be a non-degenerate braided quiver
and let $\iota: \Ac \to \gstr$ be the canonical map.

\begin{theorem}\label{accion2}
(a). There is a unique structure of braided groupoid on
$\gstr$ such that $\iota$ is a morphism of braided quivers.

\bigbreak (b). The braided groupoid $\gstr$-- with the structure
in (a)-- is universal in the following sense. If $\G$ is a braided
groupoid and $\varphi: \Ac \to \G$ is a morphism of braided
quivers then there is a unique morphism of braided groupoids
$\widehat\varphi: \gstr \to \G$ such that $\varphi =
\widehat\varphi \iota$.
 \end{theorem}

\pf (a). By Lemmas \ref{rigid-double} and \ref{rigid-paths}, there
is a structure of braided quiver on $\path (\D \Ac)$. We first
claim that it descends to the free groupoid $F(\Ac)$. Let $(u,v)
\in \path (\D \Ac ) {\,}_\tgt\hspace{-0.1cm} \times_\sou \path (\D
\Ac )$ and let $u',v'$ be elementary reductions of $u,v$
respectively. By \eqref{trenza-path1} and \eqref{trenza-path2}, we
have $u\fde v = u'\fde v$ and $u\fiz v = u\fiz v'$. Now let $x\in
\D \Ac$, $\sou(x) = \tgt(u)$. Then
$$
u\fde (xx^{-1}) = (u\fde x) \left((u\fiz x) \fde x^{-1} \right) =
(u\fde x)(u\fde x)^{-1},
$$
the first equality by \eqref{trenza-path3} and the second by
\eqref {trenzas-notac5}, since $\path (\D \Ac) = \path (\D
\Ac)^{\op}$ as braided quivers. This implies that $u\fde v \sim
u\fde v'$. Similarly $u\fiz v \sim u'\fiz v$. In conclusion,
$F(\Ac)$ is a braided quiver with structure inherited from $\path
(\D \Ac)$.

\bigbreak We next claim that this descends to the structure
groupoid $\gstr$. To see this, we observe that the kernel of the
canonical map $F(\Ac) \to \gstr$ is the subgroup bundle of
$F(\Ac)$ given by
$$\n = \left\{n_1^{\pm 1} \dots n_k^{\pm 1}: n_i
= x_iy_i(x_i\fiz y_i)^{-1} (x_i\fde y_i)^{-1},
\text{ for some } (x_i, y_i) \in \Ac{\,}_\tgt\hspace{-0.1cm} \times_\sou \Ac \right\}.
$$
It is indeed clear that $\n$ is the subgroup bundle generated by
the elements of the form $<x,y>:= xy(x\fiz y)^{-1}(x\fde y)^{-1}$,
so it remains to check that it is normal. Let $(u,x,y) \in
\Ac{\,}_\tgt\hspace{-0.1cm} \times_\sou \Ac
{\,}_\tgt\hspace{-0.1cm} \times_\sou \Ac$. Then
\begin{align*}
u<x,y>u^{-1} &=uxy(x\fiz y)^{-1}(x\fde y)^{-1}u^{-1}
= <ux,y> (ux\fde y) (ux\fiz y)(x\fiz y)^{-1}(x\fde y)^{-1}u^{-1}\\
&= <ux,y> (ux\fde y) (u\fiz (x\fde y))(x\fde y)^{-1}u^{-1} =
<ux,y> <u, x \fde y>^{-1} \in \n,
\end{align*}
where we have used \eqref{trenza-path1} and \eqref{trenza-path3}.
This implies that $\n$ is normal. We have
$n\fde y= y$, $x\fiz n= x$,
by \eqref{cond-structur1} and \eqref{cond-structur2}, and
hence $n\fiz y \equiv y^{-1}ny \mod \n$,
$x\fde n\equiv xnx^{-1}\mod \n$, if $(x,n,y) \in \Ac{\,}_\tgt\hspace{-0.1cm}
\times_\sou \n {\,}_\tgt\hspace{-0.1cm} \times_\sou \Ac$. Therefore,
the maps $\fde, \fiz$ descend to well-defined maps
$\fde, \fiz: \gstr {\,}_\tgt\hspace{-0.1cm} \times_\sou
\gstr \to \gstr$; these define a map $\sigma: \gstr {\,}_\tgt\hspace{-0.1cm}
\times_\sou \gstr \to \gstr {\,}_\tgt\hspace{-0.1cm} \times_\sou
\gstr$, and this is clearly a solution.

\bigbreak We next claim that $\gstr$ is a braided groupoid.
Indeed, \eqref{braidedgpds-1} follows by induction using
\eqref{trenza-path1}, \eqref{trenza-path2}, \eqref{trenza-path3}
and \eqref{trenza-path4}. The structure is unique since
$\iota(\Ac)$ generates $\gstr$ as a groupoid.

Finally, the proof of (b) is straightforward. \epf

\bigbreak
We are now ready to prove the main result of this paper.

\begin{definition}\label{sarandi}
A \emph{structural  pair} is a pair $(\G, \Ac)$, where $\G$ is a
braided groupoid and $\Ac$ is a representation of $(\G, \G \bowtie
\G)$ such that
\begin{flalign}
\label{sarandi1}
&\text{The image $\vert \Ac \vert$ generates the groupoid $\G$.}
\\ \label{sarandi2}
&\text{The map $\nabla: \G\to \aut (\D\Ac \ftimes \D\Ac)$ induced by
the left action is injective.}
\end{flalign}
\end{definition}

\begin{theorem}\label{finito} There is a bijective correspondence between

\begin{enumerate}
\item[(a)] Non-degenerate braided quivers.

\bigbreak
\item[(b)] Structural  pairs.
\end{enumerate}
\end{theorem}

\pf Let $(\Ac, \sigma)$ be a non-degenerate braided quiver.
By Lemma \ref{antipode2} (b) and Theorem \ref{accion2}
the structure groupoid of $(\Ac, \sigma^{-1})$ coincides with $\gstr$.
Thus we have left and right actions $\fd: \gstr \ftimes \Ac\to \Ac$,
$\fz: \Ac \ftimes \gstr \to \Ac$.
We define a left action $\fdd: \gstr \bowtie \gstr \ftimes \Ac\to \Ac$
by
\begin{equation}\label{acciont-quiver}
(g,h) \fdd x = g\fd(h \fde x),
\end{equation}
$(g,h)\in \gstr \bowtie \gstr$,  $x \in\Ac$, $\tgt(h) = \sou(x)$.
We check that $\fdd$ is a left action. By Lemma \ref{acciones-compatibles}, it is enough to verify \eqref{acciont-quiver1}:
$g \fde \left(h \fd x \right) \overset{?}= \left(g\fde h\right) \fd \left[(g\fiz h)\fde x \right]$ for $g,h,x$ as above. We know that this is true if $g,h \in \D\Ac$ by \eqref{tecnica1} applied to $\osigma$, \emph{cf.} Lemma \ref{rigid-double}. The set
$I_1 = \big\{h\in \gstr:$  \eqref{acciont-quiver1} holds for all $g \in \D\Ac$, $x \in \Ac\big\}$ is closed under multiplication. Indeed, if $h,k \in I_1$ then
\begin{align*}
g \fde \left((hk) \fd x \right) &= g \fde \left(h \fd (k \fd x) \right)
= \left(g\fde h\right) \fd \left[(g\fiz h)\fde (k\fd x) \right]
\\
&= \left(g\fde h\right) \fd \big[\left((g\fiz h)\fde k\right)\fd
\left[\left((g\fiz h)\fiz k\right) \fde x\right] \big]
\\
&= \left((g\fde h)\left((g\fiz h)\fde k\right)\right) \fd \big[\left(g\fiz hk\right) \fde x\big] = \left(g\fde hk\right) \fd \left[(g\fiz hk)\fde x \right].
\end{align*}
Hence $I_1 = \gstr$. Similarly, $I_2 = \big\{g\in \gstr:$  \eqref{acciont-quiver1} holds for all $h\in \gstr$, $x \in \Ac\big\}$ is closed under multiplication and hence equals $\gstr$. Thus $\fdd$ is a left action. The map $\iota: \Ac \to \gstr$ preserves the actions $\fde, \fiz, \fd, \fz$, hence it verifies \eqref{compcond} because of the definitions \eqref{accionbt1}, \eqref{acciont-quiver}. Thus $\Ac$ is a representation of the matched pair $(\gstr, \gstr \bowtie \gstr)$.

\bigbreak
Let $\Lambda$ be the normal subgroup bundle of $\gstr$ given by the
intersection of the kernels of the actions $\fde$ and $\fiz$ on $\D\Ac$. We define the \emph{reduced structure groupoid}
\begin{equation}\label{gpdt-quiver}
\ggstr := \gstr / \Lambda.
\end{equation}

\bigbreak
We claim that $\Lambda$ is a subgroup bundle of the normal subgroup bundle $\Gamma$ of $\gstr$ defined by \eqref{gamadef}.
Indeed, $\iota$ extends to a morphism of braided quivers $\iota: \D\Ac \to \gstr$, and it preserves the actions $\fde, \fiz, \fd, \fz$.
We have
\begin{align*}
v \in \Lambda \implies  v\fde \iota(x) = \iota(x), \quad \iota(y) \fiz v = \iota(y), \quad x,y\in \D \Ac, \quad
\sou(x) = \tgt(v) = \sou(v)  = \tgt(y).
\end{align*}
If $v\in \Lambda$ and $w\in \path (\D \Ac)$, say $w = x_1 \dots x_n$
with $x_i\in \D\Ac$, $1\le i \le n$, $\tgt(v) = \sou(w)$, then we
see by induction on $n$ that $v\fde w = w$. Similarly for the
right action $\fiz$, and thus $\Lambda\subset\Gamma$.

\bigbreak
Then $\ggstr$ is a braided groupoid by Lemma \ref{gama} and the
action \eqref{acciont-quiver} induces an action of $\ggstr \bowtie
\ggstr$ on $\Ac$ by definition of $\Lambda$. Let $\Vert \ \Vert:
\Ac \to \ggstr$ be given by $\Vert x\Vert = $ class of $\iota(x)$.
The map $\Vert \ \Vert$ also preserves the actions $\fde, \fiz, \fd, \fz$, hence it verifies \eqref{compcond}.
Thus $\Ac$ is a representation of the matched pair $(\ggstr, \ggstr \bowtie \ggstr)$.

\bigbreak We claim that $(\ggstr, \Ac)$ is a structural  pair.
Indeed, condition \eqref{sarandi1} is clear since $\iota(\Ac)$
generates $\gstr$, and condition \eqref{sarandi2} follows at once from the definition of $\Lambda$.

\bigbreak Conversely, if $(\G, \Ac)$ is a structural  pair then
$\Ac$ has a structure of braided quiver by Corollary \ref{exabrq}.

\bigbreak Let us finally check that both constructions are reciprocal.
If $\Ac$ is a braided quiver and $(a, b)\in \Ac \ftimes \Ac$,
then we compute the braiding arising from $\ggstr$ by \eqref{trenza23}:
\begin{align*}
\sigma_{\Ac, \Ac} (a, b)  &= \Bigl(\lyzdd(\Vert a\Vert) \fdd b,
\left( \lyzuu(\Vert b\Vert)^{-1} \fzz \Vert a \Vert^{-1}\right)
\fdd a \Bigr)
\\ &= \Bigl((\id \sou(a), \Vert a\Vert) \fdd b,
\left( (\Vert b\Vert^{-1}, \id \sou(b)) \fzz \Vert a \Vert^{-1}\right)
\fdd a \Bigr)
\\ &= \Bigl(\id \sou(a) \fd \left(\Vert a\Vert \fde b\right),
\left( \Vert b\Vert^{-1}\fz \Vert a \Vert^{-1}, \id \sou(a) \right)
\fdd a \Bigr)
\\ &= \Bigl(a\fde b, \left(b^{-1}\fz a^{-1}\right)\fd a \Bigr)
= \Bigl(a\fde b, \left(b^{-1}\fd a^{-1}\right)^{-1} \Bigr)
\\ &= \Bigl(a\fde b, a \fiz b \Bigr)=\sigma(a, b).
\end{align*}

\bigbreak Conversely, let $(\G, \Ac)$ be a structural pair. By Theorem \ref{accion2} there is a unique morphism of braided groupoids
$\psi: \gstr \to \G$ such that $\psi(\iota(a)) = \vert a\vert$, $a\in
\Ac$; by condition \eqref{sarandi1}, $\psi$ is surjective.
We claim that $\Ker \psi = \Lambda$. It is enough to show
that the following diagram commutes:
\begin{equation*}
\xymatrix{\ar[rd]_{\nabla} \gstr \ar[rr]^{\psi} &  & \G
\ar[dl]^{ \nabla}
\\ &  \aut \D\Ac.   }
\end{equation*}
Let $(g,x,y) \in \gstr \ftimes \D\Ac \ftimes \D\Ac$. Then
$$\psi(g)\fdd (x, y)
= (\psi(g)\fde x, (\psi(g)\fiz\vert x\vert) \fd y) = (\psi(g)\fde
x, \psi(g\fiz\iota (x)) \fd y),$$ the second equality since $\psi$
is a morphism of solutions. Hence, we are reduced to prove
$\psi(g)\fde x \overset{?}= g\fde x$, $z \fiz \psi(g) \overset{?}=
z \fiz g$, if  $(z,g,x) \in \D\Ac \ftimes\gstr \ftimes  \D\Ac$. If
$g\in \iota(\D\Ac)$, the identities hold because of the definition
of $\psi$. Then the identities hold always. \epf

\begin{obs} The notion of ``structural  pair" is a generalization and
reformulation of the notion of ``faithful bijective 7-uple" in
\cite{s}; Theorem \ref{finito} is a generalization of Soloviev's
Theorem \cite[Th. 2.7]{s}. In our formulation, Theorem 2.7 in
\cite{s} reads as follows:

``There is a bijective correspondence between non-degenerate
braided sets and  pairs $(G, A)$, where $G$ is a braided group and
$A$ is a representation of the matched pair $(G, G\bowtie G)$ such
that \eqref{sarandi1} and \eqref{sarandi2} hold."
\end{obs}

\begin{obs} If the braided quiver $\Ac$ is finite then
the reduced structure groupoid $\ggstr$ is finite by condition
\eqref{sarandi2}.
\end{obs}

\begin{obs}\label{faithful} Let us say that a non-degenerate braided quiver
is \emph{faithful} if the map $\iota$ is injective. In this case,
$\Lambda = \Gamma$ as in the proof of Theorem \ref{finito}, and
$\Vert \ \Vert$ is also injective. Thus, faithful non-degenerate
braided quivers are in bijective correspondence with pairs $(\G,
\Ac)$ where $\G$ is a braided quiver and $\Ac$ is an invariant
sub-quiver that generates $\G$.
\end{obs}

\begin{obs} The structural pair of a rack $(X, \trid)$ is
$(\Inn_{\trid} X, \phi)$ where the group $\Inn_{\trid} X$ is
braided via the adjoint representation and $\phi: X \to
\Inn_{\trid} X$ is the map $\phi(x)(y) = x\trid y$, see \cite{AG}.
\end{obs}

\begin{obs} An explicit group-theoretical description of matched pairs of groupoids can be  found in \cite[Th. 2.15]{AN}, see also \cite{AM}. It follows from
this description that a finite braided groupoid is roughly
determined by a group $D$ with two subgroups $V$ and $H$ such
that:
\begin{itemize}
\item There is a bijection $\Pc \simeq V\backslash D/H$,
\item $V$ and $H$ are isomorphic,
\item $V$ intersects trivially any conjugate of $H$.
\end{itemize}
\end{obs}
Details will appear elsewhere.

\bigbreak
\subsection{Rack bundles}

\

\bigbreak Let $p: \Lc \to \Pc$ be a loop bundle. If $\sigma: \Lc
{\,}_{p} \hspace{-0.1cm} \times_p \Lc \to \Lc {\,}_p
\hspace{-0.1cm}\times_p \Lc$ is a solution, and $\Lc_P$ is the
fiber of $P\in \Pc$ then the restriction $\sigma_P: \Lc_{P} \times
\Lc_P \to  \Lc_P \times \Lc_P$ is a solution to the
set-theoretical braid equation.  In other words, a solution $(\Lc,
\sigma)$ with $\Lc$ a loop bundle is the same as a bundle of
solutions of the set-theoretical braid equation. For instance, a
\emph{rack bundle} is a pair $(\Lc, \trid)$ where $\Lc$ is a loop
bundle and $\trid = (\trid_P)_{P\in \Pc}$, where $\trid_P$ is a
structure of rack in the fiber $\Lc_P$, $P\in \Pc$; we omit the
subscript in what follows. (See \emph{e. g.} \cite{AG} for
information on racks). This means that
\begin{align}
\label{ccc1} \phi_x: \Lc_P\to \Lc_P, \quad \phi_x (y)
&= x \trid y, \quad   \text{is a bijection for all } x\in \Lc_P, \, P\in \Pc,  \\
\label{ccc4} x \trid(y \trid z)&=(x\trid y)\trid(x\trid z)\
\qquad x,y,z \in \Lc_P, \, P \in \Pc.
\end{align}

\begin{obs}\label{rack} Let $p: \Lc \to \Pc$ be a loop bundle and
let  $\trid: \Lc {\,}_{p} \hspace{-0.1cm} \times_p \Lc \to  \Lc$
be a morphism of quivers. Let $c: \Lc {\,}_{p} \hspace{-0.1cm}
\times_p \Lc \to  \Lc {\,}_p \hspace{-0.1cm}\times_p \Lc$ be given by
\begin{equation}\label{eq:calt} c(x,y)=(x\trid y,x), \qquad (x,y)
\in \Lc {\,}_{p} \hspace{-0.1cm} \times_p \Lc.\end{equation}
Then $c$ is a non-degenerate solution if and only if $(\Lc, \trid)$
is a rack bundle.

\bigbreak Observe also that, if $c: \Ac {\,}_\tgt \hspace{-0.1cm}
\times_{\mathfrak s} \Ac \to \Ac {\,}_\tgt \hspace{-0.1cm}
\times_{\mathfrak s} \Ac$ is an isomorphism of quivers defined
by a formula analogous to \eqref{eq:calt} and with
$(\sou(x), x\trid\underline{\;\;}\,,\tgt (x))\in \aut \sou$
for any $x\in \Ac$ then $\Ac$ is a loop bundle. Indeed, given
$x\in \Ac$ there exists $y\in \Ac$ such that $x \trid y = x$
but then $\sou(x) = \tgt (x)$.
\end{obs}

If $(\Lc, \trid)$ is a rack bundle then we set
$$\aut_{\trid} \Lc := \left\{(P,x,Q): P,Q\in \Pc, \text{ and } x: \Ec_Q
\to \Ec_P \text{ is an isomorphism of racks} \right\}.$$

\begin{exa} If $\N$ is a group bundle, define
$\trid: \N {\,}_{p} \hspace{-0.1cm} \times_p \N \to  \N$ by $x\trid y := xyx^{-1}$;
then $(\N, \trid)$ is a rack bundle.
\end{exa}
\bigbreak
\subsection{The derived solution. Rack-theoretical characterization of braided quivers}

\

\bigbreak Let $\Ac$, $\widetilde \Ac$ be quivers and $\sigma:
\Ac{\,}_\tgt\hspace{-0.1cm} \times_\sou \Ac \to
\Ac{\,}_\tgt\hspace{-0.1cm} \times_\sou \Ac$,  $\widetilde \sigma:
\widetilde \Ac{\,}_\tgt\hspace{-0.1cm} \times_\sou \widetilde
\Ac\to\widetilde \Ac{\,}_\tgt\hspace{-0.1cm} \times_\sou
\widetilde \Ac$  be isomorphisms of quivers. We say that $(\Ac,
\sigma)$ and $(\widetilde \Ac, \widetilde \sigma)$ are
\emph{equivalent} if there exists a family of bijections $U^n:
\Ac^n \to \widetilde \Ac^n$ such that $U^n \sigma_{i, i+1} =
\widetilde \sigma_{i, i+1}U^n$, for all $n \ge 2$, $1\le i \le
n-1$.

\bigbreak
\begin{obs}\label{equiva}
If $(\Ac, \sigma)$ and $(\widetilde \Ac, \widetilde \sigma)$ are
equivalent and $(\Ac, \sigma)$ is a solution, then $(\widetilde
\Ac, \widetilde \sigma)$ is also a solution and the $U^n$'s
intertwine the corresponding actions of the braid group $\mathbb
B_n$.
\end{obs}

\begin{definition}\label{quiv-datum} A \emph{1-cocycle quiver datum}
is a collection $(\Ac, \Lc, \f, \mu)$ where $\Ac$ is a quiver
over $\Pc$, $\Lc$ is a rack bundle over $\Pc$, $\f: \Ac \to \aut
\Lc$ is an injective morphism of quivers over $\Pc$, and $\mu: \Lc
\to \Ac^\tgt$ is an isomorphism of bundles over $\Pc$, subject to
the cocycle condition \eqref{tresveinticinco-5} below.

For simplicity of the notation, we shall identify $\Ac$ with a
sub-quiver of $\aut \Lc$, and denote the inverse $\mu^{-1}:
\Ac^\tgt \to\Lc$  by $\mu^{-1}(x) = \x$.

To state the cocycle condition, we define first $\fiz: \Ac
{\,}_\tgt \hspace{-0.1cm}\times_\sou \aut \Lc \to \Ac$ by
\begin{align}\label{tresveinticuatro}
x \fiz y &= \mu y^{-1} (\x),
\qquad (x,y) \in \Ac {\,}_\tgt \hspace{-0.1cm}\times_\sou \aut \Lc.
\end{align}
In other words, this is the natural right action of $\aut \Lc$ on
$\Lc$ pulled back to $\Ac$ via $\mu$. Clearly $\tgt(x\fiz y) =
\tgt (y)$. By restriction we have $\fiz: \Ac \ftimes \Ac \to \Ac$.
We then define $\fde: \Ac {\,}_\tgt \hspace{-0.1cm}\times_\sou \Ac
\to \Ac$ by
\begin{align} \label{tresveinticinco}
x \fde y &= \mu \left(\overline{x \fiz y} \trid
\y \right) \fiz(x \fiz y)^{-1},
\qquad (x,y) \in \Ac {\,}_\tgt \hspace{-0.1cm}\times_\sou \Ac.
\end{align}

This map is well-defined and $\tgt(x\fde
y) = \sou (x \fiz y)$. \emph{The cocycle condition is}

\begin{equation}\label{tresveinticinco-5}
x \fde y = xy(x \fiz y)^{-1},
\qquad (x,y) \in \Ac {\,}_\tgt \hspace{-0.1cm}\times_\sou \Ac.
\end{equation}
\end{definition}

Note that \eqref{tresveinticinco-5} implies
$\sou(x\fde y) = \sou(x)$.
Note also that \eqref{tresveinticinco} is equivalent to
$$(z \fiz y^{-1}) \fde y = \mu(\z\trid \y) \fiz z^{-1},
\qquad (z,y) \in \Ac {\,}_\tgt \hspace{-0.1cm}\times_\sou \Ac.$$
This is in turn equivalent to \eqref{eq:rackasoc} below.

Now we are ready to state the main result of this subsection.

\begin{theorem} Let $\Ac$ be a quiver.
There is a bijective correspondence between
\begin{enumerate}
\item[(a)] Structures of non-degenerate braided quiver $(\Ac, \sigma)$.

\bigbreak
\item[(b)] 1-cocycle quiver data $(\Ac, \Lc, \f, \mu)$.
\end{enumerate}
\end{theorem}

\pf The proof of ``(a) $\implies$ (b)" is given by the next Lemma,
that generalizes results from \cite{s,lyz1}.

\begin{lema}\label{slyz1} Let $(\Ac, \sigma)$ be a non-degenerate braided quiver
and let $\Lc := \Ac^\tgt$ be the loop bundle as in \eqref{loop2}.
Let $\trid: \Lc {\,}_{p} \hspace{-0.1cm} \times_p \Lc \to  \Lc$ be
the morphism of bundles, resp. $\f: \Ac \to \aut \Lc$ the morphism
of quivers, defined by
\begin{align}\label{eq:rackasoc}
\x\trid \y &=  \overline{\left( (x \fiz y^{-1}) \fde y \right) \fiz x},
\\\label{eq:rackasocbis}
 \f_y(\x) &= \overline{x\fiz y^{-1}},
\end{align}
$x,y \in \Ac$, $\tgt(x) = \tgt (y)$.

\bigbreak
(a). $\f$ induces a morphism of groupoids $\f: \Gb_{\Ac} \to \aut_{\trid} \Lc$.

\bigbreak (b). $(\Lc,\trid)$ is a rack bundle. If $c$ is given by
\eqref{eq:calt}, then $c$ is a solution, called the \emph{derived
solution} of $\sigma$. The  solutions $\sigma$ and $c$ are
equivalent.
\end{lema}

\pf Note that \eqref{eq:rackasoc} and \eqref{eq:rackasocbis} are
well-defined since $\tgt(x) = \tgt (y)$. By
\eqref{cond-structur2}, $\f$ induces a morphism of groupoids $\f:
\Gb_{\Ac} \to \aut \Lc$.

\bigbreak
(a). We have to show $\f_y(\x \trid \z) \overset{?}{=} \f_y(\x) \trid \f_y(\z)$
if $\tgt(x) = \tgt (y)= \tgt (z)$. Now

\begin{align*}
\f_y(\x) \trid \f_y(\z) &= \overline{x\fiz y^{-1}} \trid \overline{z\fiz y^{-1}}
\\ &= \overline{\left( \left((x\fiz y^{-1}) \fiz (z\fiz y^{-1})^{-1}
\right) \fde (z\fiz y^{-1}) \right) \fiz (x\fiz y^{-1})}
\\ &= \overline{\left( \left((x\fiz y^{-1}) \fiz yz^{-1}(z\fde
y^{-1}) \right) \fde (z\fiz y^{-1}) \right) \fiz (x\fiz y^{-1})}
\\ &= \overline{\left( \left((x\fiz z^{-1}) \fiz (z\fde y^{-1})
\right) \fde (z\fiz y^{-1}) \right) \fiz (x\fiz y^{-1})}
\\ &= \overline{\left( \left((x\fiz z^{-1}) \fde z \right) \fiz
\left( ((x\fiz z^{-1}) \fiz z) \fde y^{-1} \right)  \right) \fiz
(x\fiz y^{-1})}
\\ &= \overline{\left((x\fiz z^{-1}) \fde z \right) \fiz  (x\fde
y^{-1})(x\fiz y^{-1})}
\\ & = \overline{\left( (x \fiz z^{-1}) \fde z \right) \fiz xy^{-1}}
\\ & = \f_y(\x \trid \z).
\end{align*}
Here the third, fourth, sixth and seventh equalities are by
\eqref{cond-structur2}, and the fifth by \eqref{cond-structur3}.

\bigbreak (b).  Let $U^n: \Ac^n \to \Lc^n$ be defined inductively
by
\begin{align}\label{defun1}
U^2(x_1, x_2) &= (\overline{x_1 \fiz x_2}, \overline{x_2}),
\\ \label{defun2} U^{n+1} &= Q_n (U^n \times \id),
\text{ where } Q_n(\overline{x_1},
\dots, \overline{x_n}, x_{n+1}) = (\overline{x_{1} \fiz x_{n+1}},
\dots, \overline{x_{n} \fiz x_{n+1}}, \overline{x_{n+1}}).
\end{align}

We claim  that
\begin{equation}\label{equivn}
c_{i,i+1}U^n
= U^n \sigma_{i, i+1}, \qquad n\ge 2, \, 1\le i \le n-1.\end{equation}
If $n=2$, then
\begin{align*}
c\,U^2(x_1, x_2) &= c (\overline{x_1 \fiz x_2}, \overline{x_2}) =
(\overline{x_1 \fiz x_2} \trid \overline{x_2}, \overline{x_1 \fiz
x_2})
\\ & = (\overline{\left( ((x_1 \fiz x_2) \fiz x_2^{-1}) \fde x_2
\right) \fiz (x_1 \fiz x_2)}, \overline{x_1 \fiz x_2})
\\ &  = (\overline{( x_1  \fde x_2 ) \fiz (x_1 \fiz x_2)}, \overline{x_1 \fiz x_2})
 \\ & = U^2( x_1  \fde x_2 , x_1 \fiz x_2) = U^2\sigma( x_1, x_2).
\end{align*}
Assume that \eqref{equivn} is true for $n$ and let $i$ be such
that $1\le i \le n-1$. Then
\begin{equation*}
c_{i,i+1}U^{n + 1} =c_{i,i+1}Q_n (U^n \times \id) = Q_n c_{i,i+1}
(U^n \times \id) =Q_n (U^n \times \id)\sigma_{i, i+1} = U^{n + 1}
\sigma_{i, i+1};\end{equation*} here the second equality follows
from part (a). So, let $i = n$. We claim  that
$c_{n,n+1}Q_n(Q^{n-1} \times \id) = Q_n (Q^{n-1} \times \id)
\sigma_{n, n+1}$; indeed
\begin{align*}
c_{n,n+1}Q_n(Q^{n-1} \times \id) & (\overline{x_1}, \dots,
\overline{x_{n-1}}, x_n, x_{n+1}) = c_{n,n+1}Q_n(\overline{x_{1}
\fiz x_{n}}, \dots, \overline{x_{n-1} \fiz x_{n}},
\overline{x_{n}}, x_{n+1})
\\ & = c_{n,n+1}(\overline{(x_{1} \fiz x_{n})\fiz x_{n+1}}, \dots,
\overline{(x_{n-1} \fiz x_{n}) \fiz x_{n+1}}, \overline{x_{n}\fiz
x_{n+1}}, \overline{x_{n+1}})
\\ & = (\overline{x_{1} \fiz x_{n} x_{n+1}}, \dots, \overline{x_{n-1}
\fiz x_{n} x_{n+1}}, \overline{( x_n  \fde x_{n+1} ) \fiz (x_n
\fiz x_{n+1})}, \overline{x_{n}\fiz x_{n+1}}),\end{align*} and
this equals
\begin{align*}
&Q_n (Q^{n-1} \times \id) \sigma_{n, n+1} (\overline{x_1}, \dots,
\overline{x_{n-1}}, x_n, x_{n+1}) = Q_n (Q^{n-1} \times \id)
(\overline{x_{1}}, \dots, \overline{x_{n-1}}, x_{n} \fde x_{n+1},
x_{n} \fiz x_{n+1})
\\ & = (\overline{x_{1} \fiz (x_{n} \fde x_{n+1}) (x_{n} \fiz
x_{n+1})}, \dots,  \overline{( x_n  \fde x_{n+1} ) \fiz (x_n \fiz
x_{n+1})}, \overline{x_{n}\fiz x_{n+1}})\end{align*} because of
\eqref{cond-structur2}. Then
\begin{align*}
c_{n,n+1}U^{n + 1} &= c_{n,n+1}Q_n (Q^{n-1} \times \id) (U^{n-1} \times \id)
\\ & = Q_n (Q^{n-1} \times \id) \sigma_{n, n+1} (U^{n-1} \times \id) =
Q_n (Q^{n-1} \times \id)  (U^{n-1} \times \id) \sigma_{n, n+1}=
U^{n + 1} \sigma_{i, i+1}.\end{align*}

Hence $\sigma$ and $c$ are equivalent; then $c$ is a solution by
Remark \ref{equiva}, and $(\Lc,\trid)$ is a rack bundle by Remark
\ref{rack}. \epf

We now prove the implication ``(b) $\implies$ (a)". This was
implicit in \cite{s}, in the set-theoretical setting. See also
\cite[Prop. 5.4 (3)]{AG}.

\bigbreak Let us consider a collection $(\Ac, \Lc, \f, \mu)$ as in
Definition \ref{quiv-datum} but without assuming the cocycle
condition \eqref{tresveinticinco-5}. Define $\fiz, \fde: \Ac
\ftimes \Ac \to \Ac$ by \eqref{tresveinticuatro},
\eqref{tresveinticinco}; assume that $\sou(x\fde y) = \sou(x)$.
Let $\sigma: \Ac \ftimes \Ac \to \Ac {\,}_\tgt
\hspace{-0.1cm}\times_\sou \Ac$ be given by \eqref{trenzas-notac}:
$\sigma(x, y) =  (x\fde y, x \fiz y)$, $(x, y) \in \Ac {\,}_\tgt
\hspace{-0.1cm} \times_\sou \Ac$.

\begin{lema}\label{slyz2}
The map $\sigma$ is a solution if and only if the
condition \eqref{tresveinticinco-5} holds.
If this happens, the solutions $\sigma$ and $c$ are equivalent, and
$\sigma$ is non-degenerate.
\end{lema}

\pf Assume that \eqref{tresveinticinco-5} holds. Let $U^n: \Ac^n
\to \Lc^n$ be defined by \eqref{defun1}, \eqref{defun2}. We claim
that \eqref{equivn} holds in the present situation too. In fact,
one can repeat the proof for $n=2$ word by word, since
\eqref{tresveinticinco} is equivalent to \eqref{eq:rackasoc}. Same
for the proof of the inductive step, $i<n$ since $\Ac \subset
\aut_{\trid} \Lc$ by hypothesis. Finally, the proof of the
inductive step, $i= n$ can also be repeated because condition
\eqref{cond-structur2} follows from \eqref{tresveinticinco-5}.
Hence $\sigma$ and $c$ are equivalent, and $\sigma$ is a solution
by Remark \ref{equiva}.

Conversely, assume that $\sigma$ is a solution. Then
\eqref{cond-structur2} holds;  writing this explicitly down,
\emph{cf.} \eqref{tresveinticuatro}, we get
\eqref{tresveinticinco-5}. \epf

This finishes the proof of the Theorem. \epf

\begin{obs} Let $(\Ac, \sigma)$, $(\widetilde \Ac, \widetilde
\sigma)$ be two non-degenerate braided quivers. A \emph{morphism
of braided quivers} is a morphism of quivers $T: \Ac\to \widetilde
\Ac$ such that $T\times T$ intertwines $\sigma$ and $\widetilde
\sigma$; that is, such that
\begin{align}\label{mor-cs1}
T(x \fde y) &= T(x) \fde T(y),
\\ \label{mor-cs2}
T(x \fiz y) &= T(x) \fiz T(y),
\end{align}
$(x, y) \in \Ac {\,}_\tgt \hspace{-0.1cm}\times_\sou \Ac$.

A morphism of quivers $T$ is a morphism of braided quivers if and
only $\overline T$ is a morphism of the associated rack bundles
and \eqref{mor-cs2} holds.
\end{obs}

\section{Linearization}

Let $\ku$ be a field. Let $\X$ be a set. We denote by $\ku \, \X$
the $\ku \,$-vector space with basis $e_X$, $X\in \X$.

\bigbreak Recall that $\Pc$ is our fixed basis of quivers and
groupoids. We consider $\ku \, \Pc$ as a  (commutative,
semisimple) algebra with multiplication $e_P e_Q= \delta_{P,Q} \, e_P$,
$P, Q\in \Pc$.

\bigbreak
\subsection{The category of bimodules}\label{bimodules}

 \

\bigbreak The tensor category $\Mpp$ of $\ku \, \Pc$-bimodules,
with tensor product $\otimes_{\ku\,\Pc}$ and unit $\uno = \ku \,\Pc$, can be
identified with the category of $\Pc \times \Pc$-graded  vector
spaces (with maps homogeneous of degree 0). If $M$ is a $\Pc
\times \Pc$-graded vector space then the grading is denoted $M =
\oplus_{P, Q\in \Pc} {\,}_{P} \hspace{-0.025cm}M_{Q}$, with
${}_{P}M_{Q} = e_P M e_Q$. If $m\in M$ then we set ${}_{P}m_{Q} =
e_P m e_Q$.
The tensor product of $M,N \in \Mpp$  is given by $M\otimes_{\ku\,\Pc}
N = \oplus_{P, Q\in \Pc} \left(\oplus_{R\in \Pc} \, {}_{P}
\hspace{-0.025cm}M_{R} \otimes {}_{R} \hspace{-0.025cm}N_{Q}\right)$.

\bigbreak If $M\in \Mpp$ then  we set $M^* = \oplus_{P, Q\in \Pc}
\, {}_P\hspace{-0.025cm}(M^*)_Q$ where
${}_P\hspace{-0.025cm}(M^*)_Q := \Hom_{\ku} ({\,}_{Q}
\hspace{-0.025cm}M_{P}, \ku)$. We shall consider the map $\ev: M^*
\otimes_{\ku\,\Pc} M \to \ku\,  \Pc$ given by
$$
\ev\left(\sum_{P, Q, R\in \Pc} {}_{P}\alpha_{R} \otimes
\,{}_{R}m_{Q}\right) = \sum_{P\in \Pc} \left(\sum_{R\in \Pc}
\langle {}_{P}\alpha_{R}, \,{}_{R}m_{P}\rangle \right) e_P,
$$
for ${}_{P}\alpha_{R} \in {}_P\hspace{-0.025cm}(M^*)_R =
\Hom_{\ku} ({\,}_{R} \hspace{-0.025cm}M_{P}, \ku)$,
${}_{R}m_{Q} \in{\,}_{R} \hspace{-0.025cm}M_{Q}$.

\bigbreak \emph{Assume that $\Pc$ is finite}. Let $\mpp$ be the
full tensor subcategory of $\Mpp$ whose objects are the $\Pc
\times \Pc$-graded vector spaces with finite-dimensional
homogeneous subspaces. Let $M\in \mpp$ and choose a basis
${\,}_{P} \hspace{-0.025cm}m^i_{Q}$ of ${\,}_{P}
\hspace{-0.025cm}M_{Q}$, $i$ running in some index set $I(P,Q)$.
Let ${\,}_{P} \hspace{-0.025cm}\alpha^i_{Q}$ be the corresponding
dual basis. We consider the map $\coev: \ku\,  \Pc \to M
\otimes_{\ku\,\Pc} M^*$ given by
$$
\coev(e_P) := \sum_{P, Q \in \Pc}\,  \sum_{i \in I(P, Q)} {\,}_{P}
\hspace{-0.025cm}m^i_{Q}\otimes \,  {\,}_{P}
\hspace{-0.025cm}\alpha^i_{Q} \in {\,}_{P} \hspace{-0.025cm}M_{Q}
\otimes \Hom_{\ku} ({\,}_{P} \hspace{-0.025cm}M_{Q}, \ku) \subset
{\,}_{P} \hspace{-0.025cm}(M^* \otimes_{\ku\,\Pc} M)_{P}. $$ Then $M^*$
is the dual of $M$ and $\mpp$ is rigid.

\begin{obs} Assume that $\Pc$ is not finite. Then we define $\mpp$ as the
full tensor subcategory  of $\Mpp$ whose objects are the $\Pc
\times \Pc$-graded vector spaces $M$ with finite-dimensional
homogeneous subspaces, and satisfying a condition of finite support:
for any $P\in \Pc$, the sets
$$
\supp {\,}_{P} \hspace{-0.025cm}(M)
= \{Q\in \Pc: {\,}_{P} \hspace{-0.025cm}M_{Q} \neq 0\}
\text{ and } \supp \, (M)_P = \{Q\in \Pc: {\,}_{Q} \hspace{-0.025cm}M_{P} \neq 0\}
$$
are  both finite. Again $\mpp$ is rigid.
\end{obs}

\bigbreak There is a tensor functor $\Lin$
from the category $\Quiv(\Pc)$ of quivers
over $\Pc$ to $\Mpp$ given by
\begin{equation}\label{lin}
\begin{CD}
\Ac @>\Lin>> \ku \, \Ac = \oplus_{P, Q\in \Pc} \, {}_{P}
\hspace{-0.025cm}(\ku \, \Ac)_{Q},\end{CD} \quad \text{where } {}_{P}
\hspace{-0.025cm}(\ku \,\Ac)_{Q} := \ku \, \Ac (P, Q).
\end{equation}
By abuse of notation, if $T: \Ac \to \Bc$ is a morphism of quivers
then we also denote by $T: \ku \, \Ac\to \ku \, \Bc$ the linear
map $\Lin T$; that is, $T(e_f) = e_{T(f)}$, $f\in \Ac$.

\bigbreak If $M \in \Mpp$ then choose a basis $\Ac (P, Q)$ of
${\,}_{Q} \hspace{-0.025cm}M_{P}$, $P, Q \in \Pc$.  The union $\Ac
:= \coprod_{P, Q \in \Pc}\Ac (P, Q)$ is a quiver and $M \simeq \ku
\, \Ac$ in $\Mpp$.

\bigbreak
\subsection{Solutions of the braid equation in the category of bimodules}

 \

\bigbreak Let $\Ac$ be a quiver over $\Pc$, let $\sigma: \Ac
{\,}_{\tgt} \hspace{-0.1cm}\times_\sou \Ac \to \Ac {\,}_{\tgt}
\hspace{-0.1cm}\times_\sou \Ac$ be an isomorphism of quivers and
let $\q: \Ac {\,}_{\tgt} \hspace{-0.1cm}\times_\sou \Ac
\to\ku^\times$ be a function. Let $\sigma^{\q}: \ku \, \Ac \otimes
\ku \, \Ac\to \ku \, \Ac \otimes \ku \, \Ac$ be given by
\begin{equation}\label{eq:bvsf}
\sigma^{\q}(e_x\otimes e_y) = \q_{x,y} \, \sigma (e_x\otimes e_y) =
\q_{x, y}\, e_{x \fde y}\otimes e_{x \fiz y}, \qquad (x, y)\in \Ac
{\,}_{\tgt} \hspace{-0.1cm}\times_\sou \Ac.
\end{equation}

\begin{lema}\label{soluciones} $\sigma^{\q}$ is a solution of
the braid equation in $\Mpp$ if and only if $(\Ac, \sigma)$ a
braided quiver and
\begin{equation}\label{cociclo-braideq}
\q_{x,y} \q_{x\fiz y, z} \q_{x\fde y, (x\fiz y) \fde z} = \q_{y,
z} \q_{x, y \fde z} \q_{x \fiz (y\fde z), y\fiz z}, \qquad (x, y,
z) \in \Ac {\,}_{\tgt} \hspace{-0.1cm}\times_\sou \Ac {\,}_{\tgt}
\hspace{-0.1cm}\times_\sou \Ac.
\end{equation}
\end{lema}

\pf Straightforward. \epf

\begin{obs} Let $(\Ac, \sigma)$ be a braided quiver.
Let $\Gamma$ be an abelian group denoted multiplicatively. A map
$\q:  \Ac {\,}_{\tgt} \hspace{-0.1cm}\times_\sou \Ac \to \Gamma$
is called a \emph{2-cocycle} if it satisfies
\eqref{cociclo-braideq}. The space of all 2-cocycles, which
contains all constant functions, is denoted $Z^2(\Ac, \Gamma)$.

\bigbreak We say that two functions $\q, \widetilde \q:  \Ac
{\,}_{\tgt} \hspace{-0.1cm}\times_\sou \Ac \to \Gamma$ are
\emph{cohomologous} if there exists a function $\ub:\Ac \to
\Gamma$ such that
\begin{equation}\label{coborde-braideq}
\q_{x,y} \ub_{x\fde y} \ub_{x\fiz y} = \widetilde\q_{x,y} \ub_{x} \ub_{x}, \qquad (x, y,
z) \in \Ac {\,}_{\tgt} \hspace{-0.1cm}\times_\sou \Ac {\,}_{\tgt}
\hspace{-0.1cm}\times_\sou \Ac.
\end{equation}
This is an equivalence relation. Furthermore, if $\q$ and
$\widetilde \q$ are cohomologous and $\q$ is a 2-cocycle then
$\widetilde \q$ is also a 2-cocycle. The quotient of $Z^2(\Ac,
\Gamma)$ by this equivalence relation is denoted $H^2(\Ac,
\Gamma)$.

\bigbreak Let $\ub:\Ac \to \ku^{\times}$ be a function and let
$\phi_{\ub}: \ku \, \Ac \to \ku \, \Ac$ be the linear map given by
$\phi_{\ub}(e_x) = \ub_x e_x$, $x\in \Ac$. Let $\q$,
$\widetilde\q$ be 2-cocycles with values in $\ku^{\times}$ related
by \eqref{coborde-braideq}.  Then the solutions $\sigma^{\q}$ and
$\sigma^{\widetilde\q}$ are intertwined by $\phi_{\ub}$.
Therefore, the computation of $H^2(\Ac, \ku^{\times})$ is
desirable.

\bigbreak The previous considerations, in the set-theoretical
case, are well-known. A brief discussion is in \cite[Lemma
5.7]{AG}; according to M. Gra\~na, P. Etingof was aware of this. A
definition of the full cohomology of set-theoretical solutions is
given in \cite{ces}.

\end{obs}

\begin{obs} Let $(M, c)$ be a solution of the braid equation in $\Mpp$
of the form $(\ku\,\Ac, \sigma^\q)$. The braided quiver $(\Ac, \sigma)$
is \emph{not} determined by $(M, c)$, see \cite{AG}.
\end{obs}

\bigbreak Assume again that $\Pc$ is finite.
Recall the definition of rigid in Remark \ref{rigid}.

\begin{lema}
If $\Ac$ is finite and $\sigma^{\q}$ is a solution of the braid
equation, then $\sigma^{\q}$ is rigid if and only if $\sigma$ is
non-degenerate.
\end{lema}

\pf Assume for simplicity that $\q = 1$ and  set $M = \ku \,\Ac$,
$c = \sigma^1$. Let $(\delta_x)_{x\in \Ac}$ be the basis of $M^*$
dual to $\Ac$. Then
$$c^{\flat}(\delta_x \otimes e_y) = \sum_{z\in \Ac: \, \sou(z)
= \tgt(y), \, \tgt(y\fde z) = \tgt (x)}  \delta_{x, y\fde z} \,
e_{y\fiz z} \otimes \delta_z, \qquad x,y\in\Ac, \, \sou(x) =
\sou(y).$$ Hence, if $y\fde \underline{\quad}$ is bijective then
$c^{\flat}(\delta_x \otimes e_y) = e_{y\fiz (y^{-1} \fde x)}
\otimes \delta_{y^{-1} \fde x}$, or equivalently
$c^{\flat}(\delta_{y\fde u} \otimes e_y) = e_{y\fiz u} \otimes
\delta_{u}$. Thus, if $\sigma$ is non-degenerate then $c^\flat$ is
an isomorphism.

For the converse, observe  that $\sigma^{\flat}(\delta_x \otimes
y) = 0$ if $x$ is not in the image of $y\fde \underline{\quad}$.
Thus if $c^{\flat}$ is an isomorphism then $y\fde
\underline{\quad}$ is surjective, and \emph{a fortiori} bijective
for all $y$. Finally if $y\fiz u = t\fiz u$ then
$c^{\flat}(\delta_{y\fde u} \otimes e_y) = c^{\flat}(\delta_{t\fde
u} \otimes e_t)$, which implies $y = t$. \epf

\bigbreak
\subsection{Face models}

 \

\bigbreak A Yang-Baxter (or star-triangular) face model is
essentially the same as a solution of the braid equation in the
category of $\ku \Pc$-bimodules, see  \cite{H} and references
therein. Thus, braided quivers equipped  with 2-cocycles with
values in $\ku^{\times}$ give rise to Yang-Baxter face models. For
completeness we restate the results of the preceding subsection in
the language of face models. We begin by a definition inspired by
\cite{H}.

\begin{definition}\label{panoplia}
Let us first say that a \emph{quiver} in a category $\C$
is a pair of arrows $\sou, \tgt: \Ag \to \Pg$ in $\C$.

A \emph{double quiver} is a
quiver in the category $\Quiv$ of all quivers. That is, in
the ``vertical and horizontal" notation, a double quiver is a pair
of morphisms of quivers $t,b: \Bg \to \Hg$, where $\Bg$ and $\Hg$
are quivers in the usual sense: $l,r: \Bg \to \Vg$, $l,r: \Hg \to
\Pc$, and $t,b$ should preserve $l,r$:
\begin{equation}\label{esquinas}
tr = rt, \qquad tl = lt, \qquad br = rb, \qquad bl = lb.
\end{equation}

In short, a double quiver is a collection of sets and maps
$$\begin{matrix} \qquad\Bg &\overset{t,b}\rightrightarrows &\Hg \qquad
\\ l,r \downdownarrows &\qquad&\downdownarrows l,r
\\ \qquad \Vg &\underset{t,b}\rightrightarrows &\Pc \qquad\end{matrix}$$
satisfying \eqref{esquinas}. By abuse of notation we shall say
that $(\Bg, \Vg, \Hg)$ is a double quiver over $\Pc$; or
alternatively that $\Bg$ is a double quiver with sides in $\Vg$
and $\Hg$; or that $\Bg$ is a double quiver with sides in $\Ac$ in
case  $\Vg = \Hg = \Ac$. An element $B$ of $\Bg$ is called an
\emph{oriented box} and depicted as a box
$$ B =
\begin{matrix} \quad t \quad \\ l \,\, \boxe \,\, r \\ \quad b\quad
\end{matrix}$$ where $t = t(B)$, $b = b(B)$, $r = r(B)$, $l = l(B)$,
and the four corners are $tl(B)$,  $tr(B)$, $bl(B)$, $br(B)$. In
this picture, we keep in mind the orientations top-to-bottom and
left-to-right. Morphisms of double quivers, or of double quivers
over $\Pc$, are defined in the standard way.

\bigbreak Let $\Pc$ be a set and $\Ac$, $\Bc$ be quivers over
$\Pc$ denoted vertically and horizontally, respectively. The
\emph{coarse double quiver} with sides in $\Ac$ and $\Bc$ is the
collection $(\Vg\boxplus \Hg, \Vg, \Hg)$ where $\Vg\boxplus \Hg$
is the set of all quadruples $\begin{pmatrix} \quad x  \quad \\  f
\quad g \\ \quad y \quad\end{pmatrix}$ with $x,y\in \Hg$, $f,g\in
\Vg$ such that
\begin{equation}\label{corners}
l(x) = t(f), \quad r(x) = t(g), \quad l(y) = b(f), \quad r(y) = b(g).
\end{equation}
Such a quadruple is called a \emph{face}. We omit the obvious
description of the arrows.

\bigbreak If $(\Bg, \Vg, \Hg)$ is a double quiver over $\Pc$ then
there are maps $\Theta: \Bg \to \Vg\boxplus \Hg$, $\Xi: \Bg \to
\Hg {\,}_r\hspace{-0.1cm} \times_l \Vg$ given by
$$
\Theta\left(\begin{matrix} \quad x \quad \\ f \,\, \boxe \,\, g \\
\quad y \quad \end{matrix} \right) = \begin{pmatrix} \quad x
\quad \\  f \quad g \\ \quad y \quad\end{pmatrix}, \qquad
\Xi\left(\begin{matrix} \quad x \quad \\ f \,\, \boxe \,\, g \\
\quad y \quad \end{matrix}\right) = (x,g), \qquad \begin{matrix}
\quad x \quad \\ f \,\, \boxe \,\, g \\ \quad y \quad \end{matrix}
\in \Bg.
$$

\bigbreak
Clearly $\Theta$ is a morphism of double quivers.

\bigbreak We shall say that $(\Bg, \Vg, \Hg)$ is \emph{thin} if
$\Theta$ is injective (any box is determined by its sides) and
$\Xi$ is surjective.

\bigbreak We shall say that $(\Bg, \Vg, \Hg)$ is \emph{vacant} if
$\Xi$ is bijective. \end{definition}

We now attach a vacant double quiver to any braided quiver.

\begin{definition}\label{panoplia-trenzada}
Let $\Ac$ be a braided quiver. The associated vacant double quiver
with sides  in $\Ac$ is the collection $\Ac\bowtie \Ac$ of faces
of the shape
$$\begin{matrix}
\quad \quad \quad x  \quad \\  x \fde g \, \boxee \,\, g
\\ \quad \qquad x \fiz g  \quad
\end{matrix}, \qquad (x,g) \in \Ac
\ftimes \Ac.$$
This is well-defined by \eqref{esquinas-braided}.
\end{definition}

We next recall the definition of face models, see for example
\cite{H} and references therein.

\begin{definition}\label{facemodels}
Let $\Ac$ be a quiver. A \emph{face model} on $\Ac$ is a pair
$(\mathfrak P, \w)$, where  $\mathfrak P$ is a thin double quiver
with sides in $\Ac$ and $\w: \mathfrak P \to \ku^{\times}$ is a
function.

\bigbreak A face model $(\mathfrak P, \w)$ induces a linear map
$c^\w: \ku \, \Ac\otimes_{\ku\,\Pc} \ku \, \Ac \to \ku \,
\Ac\otimes_{\ku\,\Pc}
\ku \, \Ac$ by
\begin{equation}\label{face-trenza}
c^\w (e_x \otimes e_g) = \sum \w\left(\begin{matrix}  \quad x
\quad \\  f\, \boxe \, g \\  \quad y \quad\end{matrix}\right)
e_f \otimes e_y,
\qquad (x,g) \in \Ac \ftimes \Ac,
\end{equation}
where the sum is over all the pairs $(f,y) \in \Ac {\,}_\tgt
\hspace{-0.1cm} \times_\sou \Ac$ such  that $\begin{matrix}  \quad
x  \quad \\  f\, \boxe \, g \\  \quad y \quad\end{matrix}\in
\mathfrak P$.

A \emph{star-triangular face model} is  face model $(\mathfrak P, \w)$ such
that $(\ku\, \Ac, c^\w)$ is a solution to the braid equation in $\Mpp$.
\end{definition}

\begin{obs} (Hayashi, \cite{Hjalg}).
Any solution to the braid equation in $\Mpp$ arises as $(\ku\,
\Ac, c^\w)$ for some face model $(\mathfrak P, \w)$, \emph{cf.}
the considerations in Subsection \ref{bimodules}.
\end{obs}

By Lemma \ref{soluciones}, any pair $(\Ac, \q)$ where
$\Ac$ is a braided quiver and $\q$ is a 2-cocycle with values in
$\ku^{\times}$ gives rise to a star-triangular face model
$(\Ac\bowtie \Ac, \w)$. Namely, set
$$
\w \left( \begin{matrix}
\quad \quad \quad x  \quad \\  x \fde g \, \boxee \,\, g
\\ \quad \quad x \fiz g  \quad
\end{matrix}\right) = \q_{x,g}, \qquad (x,g) \in \Ac
\ftimes \Ac.
$$

There exist of course star-triangular face models
that do \emph{not} arise from braided quivers.

\bigbreak
\subsection{Quasitriangular quantum groupoids}

\

\bigbreak Let $(M,c)$ be any solution of the braid equation in $\Mpp$.
By a generalization of the FRT-construction, Hayashi has shown the
existence of a coquasitriangular weak bialgebra $B(M,c)$ such that
$M$ is a $B(M,c)$-comodule and $c$ arises from the
coquasitriangular structure. If $(M,c)$ is in addition rigid, then
once can even produce a weak Hopf algebra $Hc(M,c)$ with this property \cite{H0}.
Assume that $(M, c)$ is of the form $(\ku\,\Ac, \sigma)$, where
$(\Ac, \sigma)$ is a braided quiver. We give an alternative
construction of a quasitriangular weak Hopf algebra realizing $(\ku\,\Ac, \sigma)$ as above.

\bigbreak Let $H$ be a weak Hopf algebra or quantum groupoid in
the sense of \cite{bnsz, bsz}, see also \cite{nik-v}.  The
category ${\,}_{H} \hspace{-0.025cm}\mathcal M$ of left
$H$-modules is a monoidal category with tensor product $\otimes N
:= \Delta(1) M\otimes_{\ku} N$, $M, N \in  {\,}_{H}
\hspace{-0.025cm}\mathcal M$ \cite{bsz, NTV}.

\bigbreak Let ($\Vc$, $\Hc$) be a matched pair of groupoids over
$\Pc$, with $\Vc$, $\Hc$ and $\Pc$ are finite sets.
Let $\ku(\Vc, \Hc)$ be the corresponding weak Hopf
algebra introduced in \cite{AN}, with notation of \cite{AA}. We know:

\bigbreak
\begin{itemize}
\item The category $_{\ku(\Vc, \Hc)} \m$ of left $\ku(\Vc,
\Hc)$-modules can be tensorially embedded into $\Mpp$ \cite{AA}.

\bigbreak
\item If $\Ac$ is a representation of $(\Vc, \Hc)$,
then $\ku \Ac$ is naturally a left module over $_{\ku(\Vc, \Hc)} \m$,
and $\Rep(\Vc, \Hc)$ can be tensorially embedded into
$_{\ku(\Vc, \Hc)} \m$ \cite[Prop. 5.6]{AA}.

\bigbreak
\item Let $(\lyzu, \lyzd)$ be a LYZ-pair for
($\Vc$, $\Hc$). Let $\R_{\lyzu,
\lyzd} \in \ku(\Vc, \Hc) \otimes_{\ku} \ku(\Vc, \Hc)$ be the
universal R-matrix constructed in \cite[(5.8)]{AA}.
Then $(\ku(\Vc, \Hc),
\R_{\lyzu, \lyzd})$ is a quasitriangular quantum groupoid-- in the
sense of \cite{bsz,NTV}-- by \cite[Th. 5.9]{AA}. Let
$\Ac$ be a representation of $(\Vc, \Hc)$ and let $\sigma_{\Ac,
\Ac}$ be the corresponding solution. By construction, the
linearization $\sigma_{\Ac, \Ac}: \ku\,\Ac\otimes \ku\,\Ac \to
\ku\,\Ac \otimes \ku\,\Ac$ is a solution of the braid equation in
$\Mpp$ that arises also from $\R_{\lyzu, \lyzd}$ and the induced
structure of $\ku(\Vc, \Hc)$-module on $\ku\,\Ac$.
\end{itemize}

Combining these remarks with Theorem \ref{finito}, we conclude:

\begin{proposition}
Let $(\Ac, \sigma)$ be a finite non-degenerate quiver. Then the linearization $\ku\Ac$ is a module over the weak Hopf algebra
$\ku(\ggstr, \ggstr\bowtie \ggstr)$ and $\sigma$ arises from
the universal R-matrix $\R_{\lyzuu, \lyzdd}$. \qed
\end{proposition}

\section{Appendix
\\ by Mitsuhiro Takeuchi}

\begin{proposition}\label{refaa1}  If $\Ac$ is a quiver over $\Pc$, there is a
matched pair of groupoids $(\Vc(\Ac),\Hc(\Ac))$ such that there is
a one-to-one correspondence between representations of $(\Vc,\Hc)$
on $\Ac$ and morphisms of groupoids $(\Vc,\Hc) \to
(\Vc(\Ac),\Hc(\Ac))$.
\end{proposition}

\pf The matched pair $(\Vc(\Ac),\Hc(\Ac))$ is constructed as follows.
$\Vc(\Ac)$ is the free groupoid generated by $\Ac$.
For $Q\in \Pc$, let $\X_Q$ be the set of all paths in $\Vc(\Ac)$ beginning with $Q$:
$$
\begin{CD} Q @>a_1>> @>a_2>> \dots @>a_n>>, \qquad a_i \in \Vc(\Ac).
\end{CD}
$$
For $P, Q\in \Pc$, let $\Hc(\Ac)(P, Q)$ be the set of all bijections
$f: \X_Q \to \X_P$ satisfying the following conditions:

\bigbreak
\begin{enumerate}
\item $f$ preserves the length of paths, hence we can write
$$
f(a_1, a_2,\dots, a_n) = (b_1, b_2,\dots, b_n).$$

\item $a_i$ is in $\Ac$ iff $b_i$ is in $\Ac$.

\item $a_i$ is an identity iff $b_i$ is an identity.

\item $f(a_1, \dots, a_ia_{i+1}, \dots, a_n)
= (b_1, \dots, b_ib_{i+1}, \dots, b_n)$, $1\le i < n$.

\item $f(a_1, a_2,\dots, a_{n-1}) = (b_1, b_2,\dots, b_{n-1})$,
if $n>1$.
\end{enumerate}

The groupoid $\Hc(\Ac)$ is defined by arrows $\Hc(\Ac)(P, Q)$ with obvious composition.
In particular, we have $f(a_1) = b_1$. We may write $b_1 = f\fde a_1$. This gives a left action of $\Hc(\Ac)$ on $\Ac$ by (2). There is a canonical quiver map $\Ac \to \Vc(\Ac)$. If we fix $f$ and $a_1$ as above, the correspondence
$$
(a_2,a_3,\dots, a_n) \mapsto (b_2,b_3\dots, b_n)
$$
satisfies conditions 2) -- 5). Hence it determines an element of
$\Hc(\Ac)$. If we denote this map by $f\fiz a_1$, we get a matched
pair of groupoids $(\Vc(\Ac),\Hc(\Ac))$ and a canonical
representation on $\Ac$. The one-to-one correspondence between
representations of $(\Vc,\Hc)$ on $\Ac$ and morphisms of groupoids
$(\Vc,\Hc) \to (\Vc(\Ac),\Hc(\Ac))$ is constructed easily. \epf

\bigbreak As an important application of this proposition, we have
a version of the so-called FRT construction for matched pairs of
groupoids. We consider quadruples $(\Vc,\Hc, \lyzu, \lyzd)$ where
$(\Vc,\Hc)$ is a matched pair of groupoids and $(\lyzu, \lyzd)$ is
a matched pair of rotations on it. In view of Theorem
\ref{clasifbrstr}, this is equivalent to specifying a braiding
structure on $\Rep (\Vc,\Hc)$. A morphism
$$
(\alpha, \beta): (\Vc_1,\Hc_1, \lyzu_1, \lyzd_1) \to (\Vc_2,\Hc_2,
\lyzu_2, \lyzd_2)
$$
means a morphism $(\alpha, \beta): (\Vc_1,\Hc_1) \to
(\Vc_2,\Hc_2)$ of matched pairs of groupoids such that the
following diagram commutes.

\bigbreak
$$
\begin{CD}
\Vc_1 @> {\displaystyle \lyzu_1, \lyzd_1}>>\Hc_1 \\
@A\beta AA @VV\alpha V\\
\Vc_2 @> {\displaystyle \lyzu_2, \lyzd_2}>>\Hc_2.
\end{CD}
$$

\bigbreak This is equivalent to saying that the functor
$\Res^{\beta}_{\alpha}: \Rep (\Vc_2,\Hc_2) \to \Rep (\Vc_1,\Hc_1)$
preserves the braiding structure.

\bigbreak
\begin{theorem}\label{refaa2} Let $\Ac$ be a non-degenerate
braided quiver. There is a quadruple $(\gstr, \hstr,
\lyzu_{\Ac}, \lyzd_{\Ac})$ with a canonical representation $\rho_0$
on $\Ac$ satisfying the following properties.

\begin{itemize}

\bigbreak\item[(a)] The matched pair of rotations $(\lyzu_{\Ac},
\lyzd_{\Ac})$ induces the braiding on $\Ac$.

\bigbreak\item[(b)] Let $(\Vc, \Hc, \lyzu, \lyzd)$ be a general
quadruple. Let $\rho$ be a representation of $(\Vc, \Hc)$ on $\Ac$
such that $(\lyzu, \lyzd)$ induces the braiding on $\Ac$. Then,
there is a unique morphism of quadruples $$(\alpha, \beta): (\Vc,
\Hc, \lyzu, \lyzd) \to (\gstr, \hstr, \lyzu_{\Ac}, \lyzd_{\Ac})$$
such that $(\Ac, \rho) = \Res^{\beta}_{\alpha}(\Ac, \rho_0)$.

\bigbreak\item[(c)] We have a one-to-one correspondence between
data $\Vc, \Hc, \lyzu, \lyzd, \rho$ as in (b) and morphisms of
quadruples $(\alpha, \beta): (\Vc, \Hc, \lyzu, \lyzd) \to (\gstr,
\hstr, \lyzu_{\Ac}, \lyzd_{\Ac})$.
\end{itemize}
\end{theorem}

\pf
Recall the structure groupoid $\gstr$ of $(\Ac, \sigma)$,
Definition \ref{frt}; it is a quotient of $\Vc(\Ac)$ and we write $a\sim b$ if
the clases of $a,b\in \Vc(\Ac)$ are equal in $\gstr$.

Let $\hstr$ be the maximal subgroupoid of $\Hc(\Ac)$ which is
compatible with the defining relations of $\gstr$. If $P, Q\in
\Pc$, $\hstr(P, Q)$ consists of all $f\in \Hc(\Ac)(P,Q)$ such that
for all $(a_1, a_2, \dots, a_n)$ and $(b_1, b_2, \dots, b_n)$ in
$\X_Q$, $(a_1, a_2, \dots, a_n)\sim (b_1, b_2, \dots, b_n)$ iff
$f(a_1, a_2, \dots, a_n)\sim f(b_1, b_2, \dots, b_n)$. Then
$(\gstr, \hstr)$ is a sub-matched pair of $(\Vc(\Ac),\Hc(\Ac))$.
We define a LYZ-pair $(\lyzu, \lyzd)$ for $(\gstr, \hstr)$ as
follows.

\bigbreak Let $a\in \gstr(P,Q)$ and $(b_1, b_2, \dots, b_n)\in
\X_Q$. Make the following diagrams by means of matched pairs
$(\Vc(\Ac),\gstr)$ and $(\gstr,\Vc(\Ac))$.

\bigbreak
$$
\begin{CD} @>a>> \\@Vc_1VV @VVb_1V \\ @>a_1>>
\\@Vc_2VV @VVb_2V \\ @>a_2>>
\\ \dots @. \dots\\@>a_{n-1}>> \\@Vc_nVV @VVb_nV \\ @>a_n>>
\end{CD}\, , \qquad \qquad
\begin{CD} @>d_1>> @>d_2>> @.\dots @. @>d_n>>\\@VaVV @VV \alpha_1V
@VV \alpha_2V @.@VV \alpha_{n-1}V @VV \alpha_nV \\
@>b_1>> @>b_2>> @.\dots @. @>b_n>>
\end{CD}\,
$$
where $c_i$, $a_i$, $d_i$ and $\alpha_i$ are defined inductively by
$$
c_i = a_{i-1} \fde b_i, \quad a_i = a_{i-1} \fiz b_i, \quad a_{0} = a
\quad  \text{ and }\quad
\alpha_{i-1} = d_i\fde \alpha_{i}, \qquad b_i = d_{i} \fiz \alpha_{1},
\quad \alpha_{0} = a.
$$
The maps
$$
(b_1, b_2, \dots, b_n) \mapsto (c_1, c_2, \dots, c_n)
\quad  \text{ and }\quad (d_1, d_2, \dots, d_n)
$$
belong to $\hstr(P,Q)$. Let $\lyzd(a)$ and $\lyzu(a)$ be these
maps. Then $\lyzu, \lyzd: \gstr \to \hstr$ are groupoid maps
giving rise to a LYZ-pair which induces the original braiding
$\sigma$ on $\Ac$.

\epf

\bigbreak In view of the universality (b), (c) above, one may call
the matched pair $(\gstr, \hstr)$ the \emph{FRT-construction for
matched pairs of groupoids}.

\bigbreak Relations of this construction with the construction in
the text are explained as follows. If $\Ac$ is a non-degenerate
braided quiver, the structure groupoid $\gstr$ has a braided
structure (Theorem \ref{accion2}). We have a matched pair of
groupoids $(\gstr, \gstr\bowtie \gstr)$ which has a canonical
LYZ-pair $(\lyzuu, \lyzdd)$ (Theorem \ref{lematak2}). Thus we have
a quadruple $(\gstr, \gstr\bowtie \gstr, \lyzuu, \lyzdd)$. Further
there is a canonical representation of this matched pair on $\Ac$
which induces the original braiding on $\Ac$ (Theorem
\ref{finito}). Comparing with Theorem 5.2 above, we conclude that
there is a unique morphism of groupoids
$$\alpha: \gstr\bowtie\gstr\to \hstr$$
such that $(\id, \alpha): (\gstr, \gstr\bowtie \gstr) \to (\gstr,
\hstr)$ is a morphism of matched pairs and that $\lyzu_{\Ac} =
\alpha\lyzuu$ and $\lyzd_{\Ac} = \alpha\lyzdd$.


\begin{thebibliography}{999999}

\bibitem[AA]{AA}{\sc M. Aguiar} and {\sc N. Andruskiewitsch},
\emph{Representations of matched pairs of groupoids and
applications to weak Hopf algebras},  {\tt math.QA/0402118}
(2004), Contemp. Math., to appear.

\bibitem[AG]{AG}
{\sc N. Andruskiewitsch} and {\sc M. Gra\~na}, \emph{From racks to
pointed Hopf algebras}, Adv. Math. {\bf 178} (2003), 177-243;
\texttt{math.QA/0202084}.


\bibitem[AM]{AM}
{\sc N. Andruskiewitsch} and {\sc M. Mombelli},
\emph{Examples of weak Hopf algebras arising from vacant double groupoids}, {\tt math.QA/0405374} (2004).

\bibitem[AN]{AN} {\sc N. Andruskiewitsch} and {\sc S. Natale},
\emph{Double categories and quantum groupoids},  {\tt
math.QA/0308228} (2003).


\bibitem[BNS]{bnsz} {\sc G. B\" ohm}, {\sc F. Nill} and {\sc K. Szlach\' anyi},
\emph{Weak Hopf algebras I. Integral theory and $C^*$-structure},
J. Algebra {\bf 221} (1999),  385--438.

\bibitem[BS]{bsz} {\sc G. B\" ohm} and {\sc K. Szlach\' anyi},
\emph{A coassociative  $C^*$-quantum group with nonintegral dimensions},
Lett. in Math. Phys. {\bf 35} (1996),  437--456.


\bibitem[CES]{ces} {\sc J. Scott Carter, M. Elhamdadi} and {\sc m. Saito},
\emph{Homology theory for the Set-theoretical Yang--Baxter and
Knot invariants from generalizations of quandles\/},
\texttt{math.GT/0206255}.



\bibitem[D]{Dr} {\sc V. G. Drinfeld},
\emph{On some unsolved problems in quantum group theory\/},
in Quantum groups (Leningrad, 1990), Lecture Notes in Math. \textbf{1510},
Springer, Berlin, (1992), 1--8.

\bibitem[EGS]{egs} {\sc P. Etingof, R. Guralnik} and {\sc A. Soloviev},
\emph{Indecomposable set-theoretical solutions to the Quantum
Yang--Baxter Equation on a set with prime number of elements\/},
J. Algebra \textbf{242} 2 (2001), 709--719.

\bibitem[ESS]{ess} {\sc P. Etingof, T. Schedler} and {\sc A. Soloviev},
\emph{Set-theoretical solutions to the Quantum Yang--Baxter
Equation\/}, Duke Math. J. \textbf{100} (1999), 169--209.


\bibitem[FJK]{fjk} {\sc R. Fenn, M. Jordan-Santana} and {\sc L. Kauffman},
\emph{Biracks and virtual links}, preprint available at
\texttt{http://www.maths.sussex.ac.uk////Staff/RAF/Maths/}


\bibitem[FRS]{frs} {\sc R. Fenn, C. Rourke} and {\sc B. Sanderson},
\emph{An introduction to species and the rack space}, in Bozhüyük,
Topics in knot theory, Mehmet Emin (ed.),  Dordrecht: Kluwer
Academic Publishers. NATO ASI Ser., Ser. C, Math. Phys. Sci. 399,
33-55 (1993).

\bibitem[H1]{H0} {\sc T. Hayashi},
\emph{Quantum groups and quantum semigroups\/}, J. Algebra
\textbf{204} (1998), 225--254.

\bibitem[H2]{H} {\sc T. Hayashi},
\emph{A brief introduction to face algebras\/}, in  ``New trends
in Hopf Algebra Theory"; Contemp. Math. \textbf{267} (2000),
161--176.

\bibitem[H3]{Hjalg} {\sc T. Hayashi},
\emph{Coribbon Hopf (face) algebras generated by lattice
models\/}, J. Algebra \textbf{233} (2000), 614--641.

\bibitem[LYZ1]{lyz1} {\sc Jiang-Hua Lu, Min Yan} and {\sc Yong-Chang Zhu},
\emph{On Set-theoretical Yang--Baxter equation\/}, Duke Math. J.
\textbf{104} (2000), 1--18.

\bibitem[LYZ2]{lyz3} {\sc Jiang-Hua Lu, Min Yan} and {\sc Yong-Chang Zhu},
\emph{Quasi-triangular structures on Hopf algebras with positive
bases\/}, in ``New trends in Hopf Algebra Theory"; Contemp. Math.
\textbf{267} (2000), 339--356.

\bibitem[M]{mk1} {\sc K. Mackenzie},  \emph{Double Lie algebroides
and Second-order Geometry, I}, Adv. Math. {\bf 94} (1992), pp.
180--239.

\bibitem[NTV]{NTV} {\sc D. Nikshych, V. Turaev,} and {\sc L. Vainerman},
\emph{Invariants of knots and 3-manifolds and quantum groupoids},
Proceedings of the Pacific Institute for the Mathematical Sciences
Workshop "Invariants of Three-Manifolds"
(Calgary, AB, 1999).
Topology Appl. \textbf{127} (2003),  no. 1-2, 91--123.

\bibitem[NV]{nik-v} {\sc D. Nikshych} and {\sc L. Vainerman},
\emph{Finite quantum groupoids and their applications}, in
``Recent developments in Hopf algebra Theory", MSRI Publications
43 (2002), 211--262, Cambridge Univ. Press.

\bibitem[SW]{sw} {\sc D. Silver} and \textsc{S. Williams},
\emph{A generalized Burau representation for string links},
Pacific J. Math. \textbf{197} (2001), 241–255.


\bibitem[S]{s} {\sc A. Soloviev},
\emph{Non-unitary set-theoretical solutions to the quantum
Yang--Baxter equation\/}, Math. Res. Lett. \textbf{7} (2000), no.
5-6, pp. 577--596.

\bibitem[T]{tak}  {\sc M. Takeuchi},
\emph{Survey on matched pairs of groups. An elementary approach to
the ESS-LYZ theory}, Banach Center Publ. \textbf{61} (2003),
305--331.

\bibitem[W]{w}  {\sc M. Wada},
\emph{Twisted Alexander polynomial for finitely presentable
groups}, Topology \textbf{33}   (1994), pp. 241--256.


\bibitem[WX]{wx}  {\sc A. Weinstein} and {\sc P. Xu},
\emph{Classical solutions of the quantum Yang-Baxter equation},
Commun. Math. Phys. \textbf{148}   (1992), pp. 309--343.


\end{thebibliography}
\end{document}